\documentclass[11pt,a4paper,reqno]{amsart}
 \usepackage[utf8]{inputenc}
 \usepackage[english]{babel}
 \usepackage[normalem]{ulem}
\usepackage{hyperref}
\usepackage{color,latexsym,amssymb,amsmath,amsthm,amsfonts,enumerate,pdfsync,float,lmodern,multicol,paralist,graphics,graphicx,epsfig,tikz,dsfont,mathrsfs,esint,subfig,enumitem,pifont,cleveref,textcomp}
\usepackage{tikz-cd} 
\usepackage{mathtools}
\usepackage{euscript}
\usepackage[left=2cm,right=2cm,top=2.25cm,bottom=2.25cm]{geometry}
\usetikzlibrary{patterns}

\newcommand{\sfL}{\mathsf{L}}
\newcommand{\dd}{{\, \mathrm d}}

\def\eps{{\varepsilon}}
\def\D{{\mathcal{D}}}
\def\I{{\mathcal{I}}}
\def\R{{\mathbb{R}}}
\def\Z{{\mathbb{Z}}}
\def\N{{\mathbb{N}}}
\def\E{{\mathbb{E}}}
\def\ds{\displaystyle}

\newcommand{\opdiff}[1]{\omega[{#1}]}
\newcommand{\opdb}[2]{\D_{{#2}}[{#1}]}

 \newtheorem{theorem}{Theorem}[section]
 \newtheorem{prop}[theorem]{Proposition}
 \newtheorem{lemma}[theorem]{Lemma}

\title[Travelling waves in nonlocal reaction-dispersion equations]{Travelling waves in nonlocal reaction-dispersion equations with diffuse Lévy measures}

\author{\'Emeric Bouin}
\address[E. Bouin]{CEREMADE - Université Paris-Dauphine, PSL Research University, UMR CNRS 7534, Place du Mar\'echal de Lattre de Tassigny, 75775 Paris Cedex 16, France.}
\email{bouin@ceremade.dauphine.fr}

\author{Jérôme Coville}
\address[J. Coville]{UR 546 Biostatistique et Processus Spatiaux, INRAE, Domaine St Paul Site Agroparc, F-84000 Avignon, France.}
\address{ICJ UMR 5208 - Universite Claude Bernard Lyon 1, Campus de la Doua, 69622 Villeurbanne, France.}
\email{jerome.coville@inrae.fr}

\begin{document}
\begin{abstract}
In this paper, we focus on the existence of travelling fronts for nonlocal reaction-dispersion models. Our aim is to provide a unified existence proof in a very broad framework using simple real analysis tools. In particular, we review the existing literature and fill in the gaps. It appears that the central case is the bistable nonlinearity, which we then extend to other classical nonlinearities.
\end{abstract}

 \maketitle
\tableofcontents
 \section{Introduction and main results}
Let $\mu$ be a nonnegative, symmetric and diffuse Radon measure on $\R^{\star}:=\R\setminus\{0\}$ and $f$ a given function. In this article, we focus on equations of the following form:
\begin{equation}\label{eq:main}
U_t =\opdb{U}{\mu} +f(U), 
 \end{equation}
where $\D_{\mu}$ is the following operator,
\begin{equation*}
 \opdb{U}{\mu}:=\textrm{P.V.}\left(\int_{\R} [U(\cdot+z) - U(\cdot)] \, \dd \mu(z)\right)=\lim_{\delta\to 0}\int_{\{|z|\ge \delta\}} [U(\cdot+z) - U(\cdot)] \, \dd \mu(z),
\end{equation*}
with domain $\text{Dom}(\D_{\mu})$  given by 
$$\text{Dom}(\D_{\mu}):=\{u :\R\to\R\,|\, \forall x \in \R\,,  \opdb{u}{\mu}(x) \text{ exists and is finite}\}.$$
Throughout the paper, without further notice, the non-zero measure $\mu$ is assumed to satisfy
\begin{equation*}
    \int_{\R}\min(1,z^2)\,\dd \mu(z) <+\infty, \qquad\text{and} \qquad \mu(\{x\})=0 \quad \text{ for all }\quad x\in\R^{\star},
\end{equation*}
whereas the function $f$ will satisfy 
\begin{equation*}
    f\in \mathscr{C}^1(\R), \quad f(0)=f(1)=0, \quad f'(1)<0,
\end{equation*}
and be of one of the following types, commonly encountered in the literature:
\begin{itemize}
\item[$\star$] Bistable:
\begin{equation}\label{eq:bistable}
    \exists \,\theta\in(0,1), \text{ such that } f'(0)<0, \; f|_{(0,\theta)}<0, \; f(\theta)=0, \; f|_{(\theta,1)}>0, \;  f'(\theta)>0.\tag{Bi}
\end{equation}

\item[$\star$] Ignition: 
\begin{equation}\label{eq:ignition}
    \exists \, \theta\in(0,1), \text{ such that } \; f|_{[0,\theta]}\equiv 0,\qquad f|_{(\theta,1)}>0\qquad \text{and }
f(1)=0.\tag{Ig}
\end{equation}
\item[$\star$] Monostable: 
\begin{equation}\label{eq:monostable}
   f>0\text{ in }(0,1).\tag{Mo}
\end{equation}
\end{itemize} 
A travelling wave solution to \eqref{eq:main} is a particular type of entire solution of the form $U(t,x):= u(x-ct)$ where the wave speed is usually nonzero  ($c\in\R^{\star}$), $u$ is bounded ($u \in \mathsf{L}^{\infty}(\R)$) and monotone decreasing. When the speed $c$ equals $0$, we refer to the notion of standing wave.
The pair  $(c,u)$  satisfies
\begin{align}\label{eq:TW}
\begin{dcases}
\opdb{u}{\mu} +cu' + f(u) =0,\\[5pt]
\lim_{x \to - \infty} u(x)= 1, \qquad \lim_{x \to + \infty} u(x) = 0.
\end{dcases}\tag{TW}
\end{align}
 
This type of solution is a key tool to analyse propagation phenomena in various contexts ranging from phase transitions in materials to biological invasions, nerve propagation, and epidemiology. See, for example, \cite{DeMasi1994,DeMasi1994a,DeMasi1995,Ermentrout1993,Fife1996,Schumacher1980,Medlock2003,Clark1998,Clark1998a,Murray1993,Okubo2002,Kawasaki1997} and references therein.

For these three types of nonlinearities, the existence of travelling fronts has been obtained for various diffuse measures $\mu$ such that  
\begin{equation*}
     d\mu(z)=J(z)\dd z, \quad \text{ where either } \quad J\propto |\cdot|^{-(1+2s)} \text{ or } J\in \mathsf{L}^{1}(\R).
\end{equation*}
These two particular types of measures correspond to the two main classes of nonlocal operators widely encountered in the literature. They correspond respectively to the fractional Laplacian and the convolution operator. As mentioned above, the principal aim of this work is to unify and extend propagation results that have been obtained over the last two decades.     

Because the literature is
extensive and the hypotheses
ensuring the existence of
fronts vary substantially from one nonlinearity to another,  we choose to analyse each type of nonlinearity separately, beginning with the bistable situation, which will appear to be the core of the paper. 

When $f$ is bistable, existence of travelling front solutions has been proved when the measure $\mu$ is a \textit{bounded} Borel measure, thanks to the works of Alberti and Bellettini \cite{Alberti1998}, Bates, Fife, Ren and Wang \cite{Bates1997}, Chen \cite{Chen1997}, Fang and Zhao \cite{Fang2015} and Yagisita \cite{Yagisita2009a}. Similarly, when the measure $\mu(z) \sim |z|^{-(1+2s)} \dd z$, comparable existence results have also been obtained. We refer to Chmaj \cite{Chmaj2013}, Achleitner and Kuehn \cite{Achleitner2015}, Cabr\'e and Sire \cite{Cabre2015}, Gui and Zhao \cite{Gui2015a} and Palatucci, Savin, and Valdinoci \cite{Palatucci2013}.  Although it is not the topic of this paper, we mention that Alberti and Bellettini \cite{Alberti1998} discussed the existence of stationary solutions to \eqref{eq:main} when the nonlinearity is well balanced (\textit{i.e.} $\int_{0}^1 f(s)\dd s=0$), and gave an almost complete answer for generic absolutely continuous measures $\mu$, i.e. $\mu(z)=J(z)\,\dd z$. They proved the existence of a monotone standing profile  connecting  $0$ and $1$ when $J$ satisfies a first moment condition, that is
$$\int_{\R}\min\{z^2,|z|\}\,J(z)\dd z<+\infty.$$
It seemed to us that the methodology used there goes beyond the case of absolutely continuous measures $\mu$ and may be adapted to a general diffuse measure, extending the scope of the validity of this existence result.  

In view of the literature described above, we are led to the conclusion that for a generic measure $\mu$, a positive front (possibly discontinuous) with a unique speed should exist. We show that this is actually true when the measure $\mu$ is diffuse and $f$ is not well balanced. More precisely, we prove 

\begin{theorem}\label{thm:bi}
Assume that $f$ is of type \eqref{eq:bistable} and that it is not well balanced, with
$$
\int_{0}^1f(s)\, \dd s > 0.
$$
Then, there exists $(c,u) \in \R_{+} \times \mathsf{L}^{\infty}(\R)$ satisfying \eqref{eq:TW}.  In addition, $c$ is unique and the profile $u$ is monotone non increasing. If the measure is unbounded \footnote{An unbounded measure $\mu$ is a positive measure such that for  any $r>0, \mu((-r,r))=+\infty$.}, then  $c>0$  and the profile $u$ is  at least in $\mathscr{C}_{\mathrm{loc}}^{1,\alpha}(\R)$ for $\alpha\in(0,\frac{1}{2}
)$ and unique up to translation, i.e. if $v$ is another profile, then $u\equiv v(\cdot+h)$ for some $h\in\R$.
\end{theorem}

\textit{Stricto sensu}, not well balanced means that $\int_0^1 f(s) \, \dd s \neq 0$. However, one may assume positivity of $\int_0^1 f(s) \, \dd s $ \textit{w.l.o.g.}, since, if not, the bistable function $g(\cdot) = -f(1-\cdot)$ satisfies this condition and  one may recover the original problem by considering $(1-u(- \, \cdot), -c)$ where $u$ is a front defined for the nonlinearity $g$.

The above results can be seen as a partial extension of Yagisita's existence result on bounded Borel measures \cite{Yagisita2009a} to more generic Lévy measures.  
In this paper, bistable nonlinearities turn out to be very central since we carry out a full analysis in this case and then deduce results for other nonlinearities. 

Let us now address the question of the existence of front solutions to \eqref{eq:main} when the nonlinearity is of ignition type \eqref{eq:ignition}. For such $f$, travelling fronts have first been constructed for a general bounded Borel measure $\mu$ having a first moment in papers by the second author \cite{Coville2003a,Coville2006,Coville2007d} and then for the fractional Laplacian, simultaneously,  by Chmaj \cite{Chmaj2013} and Mellet, Roquejoffre and Sire \cite{Mellet2014}. 

\begin{figure}[h!]
\centering
\begin{tikzpicture}[scale=0.8]
%
\draw[-,line width= 2, color=red](0,-0.02)--(0,2.7) ;
\draw[line width=2, color=red](5,0)--(5,2.7); 
\draw[->,line width=2](-0.5,0)--(10.5,0) node[below]{$s$};

\fill[pattern=north west lines, pattern color=purple, opacity=0.25] 
(0,0)--(0,2.7) -- (5,2.7)--(5,0)-- cycle;


\fill[pattern=north west lines,pattern color=black!10!blue, opacity=0.25] 
(5,0)--(5,2.7)--(10,2.7)--(10,0)-- cycle;

\node at (0,0)[below]{$s=0$};
\node at (5,0)[below]{$s=\frac{1}{2}$};
\node at (7.5,0.5)[line width=6,blue]{ \ding{203}};
\node at (7.5,1.55)[blue]{ \textbf{Existence of }};
\node at (7.5,1) [blue]{ \textbf{front solutions}};
\node at (2.5,0.5)[line width=6,purple]{\ding{202}};
\node at (2.5,1.35)[red]{\Large \textbf{{\bf $x_\lambda(t) \asymp_\lambda  t^{\frac{1}{2s}}$}}};
\node at (4.7,1.35)[red, rotate=90]{\scriptsize\textbf{{\bf $x_\lambda(t) \asymp_\lambda  t \ln(t)$}}};

\end{tikzpicture}
\caption{Spreading in \eqref{eq:main} with a nonlinearity of type \eqref{eq:ignition} for some specific measure $\mu$, typically $\mu\sim\frac{\dd z}{|z|^{1+2s}}$.}
\label{fig:diagramig} 
\end{figure}
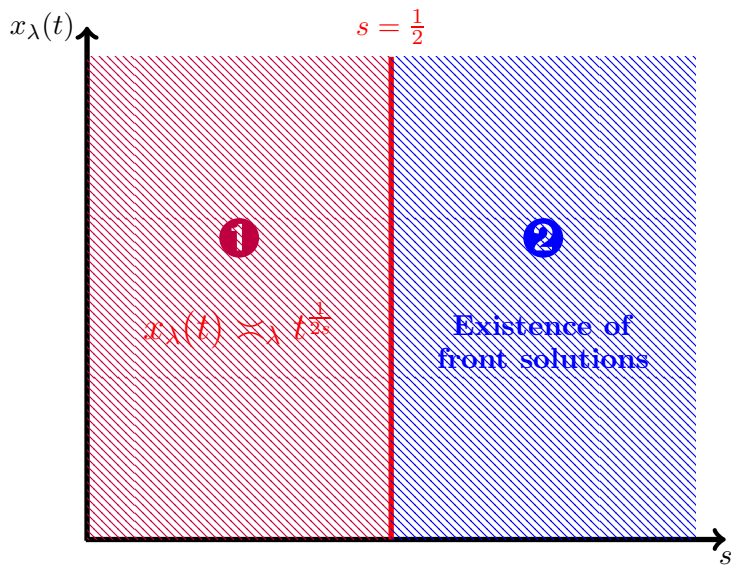

We point out that the picture is a bit different here than in the bistable framework, since in some situations, fronts structurally cannot exist. Indeed, nonexistence has been proved by Gui and Huan \cite{Gui2015} for the fractional Laplacian when $s \le \frac{1}{2}$. In such a situation, the Cauchy problem starting with Heaviside-type initial data accelerates. 
Since acceleration phenomena are not the main purpose of this discussion, we refer to \cite{Coville2021a,Zhang2023,Bouin2021ign,bouin2026propagation} for details on the spreading rates. \Cref{fig:diagramig} displays a  typical picture of spreading in such a case.

Nevertheless, and in particular from \cite{Bouin2021ign}, the existence of a front with a unique speed should exist as long as a first-moment-like condition is imposed on the measure $\mu$. With such an assumption, we prove the following theorem.

\begin{theorem}\label{thm:ign}
Assume that $\mu$  satisfies further
\begin{equation*}
    \int_{\R}\min(|z|,z^2)\,\dd \mu(z) <+\infty.
\end{equation*}
and $f$ is of type \eqref{eq:ignition}. Then there exists $(c,u)\in \R_{+}^{\star} \times \left( \mathsf{L}^{\infty}(\R)\cap \mathscr{C}^1(\R)\right)$ that satisfies \eqref{eq:TW}. In addition, $c$ is positive and unique and the profile $u$ is monotone non increasing and unique up to translation.
\end{theorem}

Note that the moment condition is identical to the one introduced by Alberti and Bellettini \cite{Alberti1998} to analyse the existence when the nonlinearity is bistable. Note also that contrary to the bistable case, the front is always at least $\mathscr{C}^1(\R)$ in this situation. This is due to the fact that for ignition type nonlinearities the speed is necessarily a strictly positive quantity when the measure $\mu$ is symmetric, regardless of whether the measure considered is bounded or unbounded.

Last, let us focus on the situation when $f$ is of type \eqref{eq:monostable}. In this situation, the existence of travelling front solutions to \eqref{eq:main} is strongly related to the interplay between the behaviour of $\mu$ and $f$. 
For instance, when $f$ is a Fisher-KPP nonlinearity, it is known that when the measure $\mu$ is absolutely continuous with a positive continuous exponentially bounded kernel $J$, that is
\begin{equation*}
\exists \,\eta \in \R_{+}^{\star}, \qquad \int_{\R} J(z) e^{\eta\vert z \vert} \, \dd z < \infty,
\end{equation*}
travelling waves do exist \cite{Schumacher1980,Weinberger1982,Carr2004,Coville2007a,Coville2008a,Lutscher2005,Yagisita2009}. 
In fact, in such a situation, the existence of fronts holds as well for more general monostable nonlinearities, see  \cite{Coville2007a,Coville2008a,Coville2007d}.  
On the other hand, when the kernel $J$ possesses heavy tails, in the sense
\begin{equation*}
\forall \, \eta \in \R_{+}^{\star},  \qquad \lim_{\vert z \vert \to \infty} J(z) e^{\eta \vert z \vert} \, = \infty,
\end{equation*}
travelling waves do not exist and the solutions exhibit an acceleration phenomenon, see \cite{Medlock2003,Yagisita2009,Garnier2011,Bouin2018,Cabre2013}. 

\begin{figure}[H]\centering
\begin{tikzpicture}[scale=0.8]
\fill[pattern=north west lines,pattern color=black!50!green, opacity=0.5] 
plot [domain=3.429:10] (\x, {(2*(\x/6))/(2*(\x/6)-1)})
--(10,8)--(3.429,8)
-- cycle;
\fill[pattern=north west lines, pattern color=blue, opacity=0.5] 
(0.5,1)--(0.5,8)--(3.429,8)
-- plot [domain=3.429:10] (\x, {(2*(\x/6))/(2*(\x/6)-1)})
-- (10,1)
-- (7,1) -- (7,1.4) -- (4.0,1.4) -- (4.0,1) 
-- cycle;
%
%
\draw[->,line width=1.2, color=black](9.5,1)--(10.5,1) node[below]{$s$};
\draw[->,line width=1.2, color=black](0.5,1)--(0.5,8.5) node[right]{$\beta$};
\draw[line width=1, dashed](3,4.2)--(3,8.5) node[right]{$s=\frac12$};
\draw[line width=1, dashed](3,1)--(3,3) ;
%
%
\draw [domain=3.429:10,line width=3, samples=150, black!50!green] plot (\x, {(2*(\x/6))/(2*(\x/6)-1)}) node[anchor=south east]{$\qquad \beta = \frac{2s}{2s-1}$};
\draw [line width=3,orange] (7,1)--(10,1) node[anchor=north east]{$\beta = 1$};
\draw [line width=3,orange] (0.5,1)--(4,1) ;
\node at (5.5,1)[line width=6,orange]{$x_\lambda(t) \asymp e^{\rho t}$};
\node at (2.25,3.5)[line width=6,blue]{$x_\lambda(t) \asymp t^{\frac{\beta}{2s(\beta-1)}}$};
\node at (6.6,5.2)[line width=15,black!50!green]{ \textbf{Existence of }};
\node at (6.6,4.7)[line width=15,black!50!green]{ \textbf{ front-like solutions}};
\end{tikzpicture}
\caption{(Rates of) propagation when $\mu(z)=J(z)\,\dd z$ with $J\in \mathsf{L}^{1}(\R)$ or $J \propto |\cdot|^{-(1+2s)}$.}
\label{fig:diagram}
\end{figure}
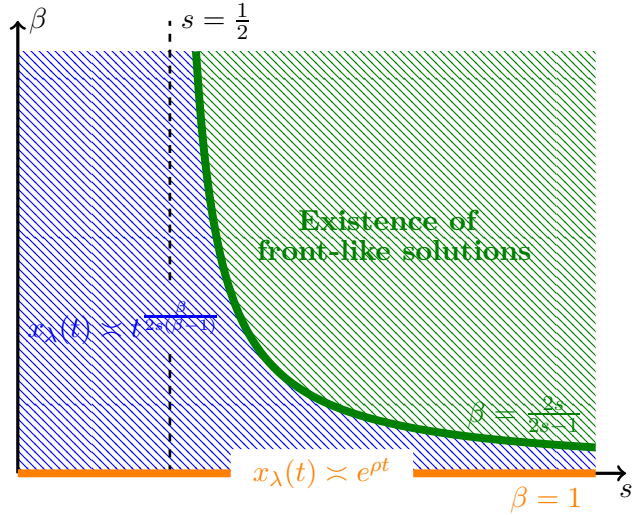

When a weak Allee effect is introduced, the study of propagation becomes more subtle. For instance, when $f$ is of the form $f(u):=u^{\beta}(1-u)$ with $\beta>1$, after a collection of works \cite{Gui2015a,Bouin2021b,Bouin2026,Alfaro2017,Coville2021a,Zhang2023}, it turns out that existence or nonexistence of travelling waves depends on whether $(2s-1)(\beta-1)\ge 1$ or not. We summarise this in \Cref{fig:diagram}.

From earlier works (see \textit{e.g.} \cite{Coville2008a}), it turns out that a common feature for proving existence of travelling fronts is the existence of special super-solutions. We prove the following result. 

\begin{theorem} \label{thm:mono}
Assume that $\mu$  satisfies further
\begin{equation*}
    \int_{\R}\min(|z|,z^2)\,\dd \mu(z) <+\infty.
\end{equation*} and that $f$ is of type \eqref{eq:monostable}. Assume further that there exists $(\kappa,w)\in \R_{+}^{\star} \times \left(\mathscr{C}^2(\R)\cap \sfL^{\infty}(\R)\cap \mathrm{Dom}(\D_\mu)\right)$ such that $w$ is positive and,
\begin{equation}\label{hyp:fmono}
\begin{dcases}
\opdb{w}{\mu}+\kappa w'+f(w) <0,\\[5pt]
\lim_{x \to - \infty} w(x) \geq 1, \qquad \inf_{x\in\R} w(x)<1.
\end{dcases}
\end{equation} 
There exists $0<c^\star\le \kappa$ such that, for $c\ge c^\star$, there exists $u\in \mathsf{L}^{\infty}(\R)\cap \mathscr{C}^1(\R)$  positive and monotone decreasing, such that $(c,u)$ satisfies \eqref{eq:TW}.
\end{theorem}

The super-solution $w$ is in particular a way to encode that acceleration cannot occur in a reasonable Cauchy problem associated to \eqref{eq:main}. 
\subsection{Strategy of proof and further comments}
Our main strategy to construct a front, whatever the nonlinearity is, is to use an approximation procedure in the spirit of the approach of Chmaj \cite{Chmaj2013}. Since existence results have been obtained for bounded Borel measures (see \cite{Yagisita2009a} for \eqref{eq:bistable}, \cite{Coville2007d} for \eqref{eq:ignition} and \cite{Yagisita2009} for \eqref{eq:monostable}), we start with $(c_\eps,u_\eps)_{\eps >0}$, solution to 
\begin{align}\label{eq:ueps}
\begin{dcases}
-c_\eps u^{\prime}_{\eps}= \opdb{u_\eps}{\mu_\eps} + f(u_{\eps}) \quad \textit{ a.e. },\\[5pt]
\lim_{x \to -\infty} u_\eps(x) = 1,\qquad \lim_{x \to +\infty} u_\eps(x) = 0,
\end{dcases}
\end{align} 
The whole question is to prove rigorously that, up to subsequences, $(c_\eps,u_\eps)_{\eps >0}$ converges towards a solution to \eqref{eq:TW}. To study this singular limit, we prove uniform bounds on $\left( c_{\eps} \right)_{\eps >0}$, as well as uniform regularity estimates on $\left(u_{\eps}\right)_{\eps >0}$ that make the family $(c_{\eps},u_{\eps})_{\eps >0}$ be pre-compact in some well chosen Banach space. Since for all $\eps > 0$, $u_\eps$  is uniformly bounded and monotone, the family $\left( u_\eps\right)_{\eps >0}$ is pre-compact for the pointwise topology thanks to Helly's theorem. However, understanding sub-sequential limits of the family  $\left( \opdb{u_{\eps}}{\mu_\eps} \right)_{\eps > 0}$ is quite delicate in general. In particular, although pointwise convergence of sub-families of $\left( u_{\eps} \right)_{\eps >0}$ would be sufficient to address this question when $\mu$ is a bounded Borel measure, it is far from being sufficient when dealing with a generic diffuse Lévy measure. Overcoming this issue is the core of the proofs presented in \Cref{sec:bi}, \Cref{sec:ign} and \Cref{sec:mono}, and their structure will be explained there in full detail.
\medskip

Before going to the proofs let us make some more comments. First, let us recall that the operator $\D_\mu$ can be seen as the infinitesimal generator of a Feller semigroup. Indeed, following \cite{Jacob2002,Jacob-Schilling2001}, the symbol $$\ds{\psi(\xi):=\int_{\R}\left(e^{iz\xi}-1-iz\xi\mathds{1}_{\{|z|\le 1\}}(z)\right)\dd \mu(z)},$$ is the Lévy-Khintchine representation of a pure jump Lévy process  $(X_t)_{t>0}$ whose characteristic function is $\varphi_{X}(\xi)(t)=\E[e^{i\xi X(t)}]=e^{t\psi(\xi)}$ .  
When $\mu$ is symmetric, then $\ds{\psi(\xi)=\int_{\R}(\cos(z\xi)-1)\dd\mu(z)}$ and 
$$\opdb{\varphi}{\mu} =\mathscr{F}^{-1} \left[\psi\mathscr{F}[\varphi]\right](x)$$
for sufficiently regular functions $\varphi$, see for example \cite{Jacob-Schilling2001, schilling2001}. As shown for example in \cite{Farkas2001}, the symbol $\psi$ completely characterises the infinitesimal behaviour of the associated Feller process and determines the regularity properties of the functions belonging to the domain $\text{Dom}(\D_{\mu})$. 

\medskip

Next, let us comment on the symmetry assumption we imposed on the measure $\mu$. This condition  ensures that the sign of the speed is uniquely related to properties of the nonlinearity, in particular depends only on  the sign of $\int_{0}^1f(s)\dd s$. Besides, it ensures that $\D_\mu$ is self adjoint in the adequate functional space (\textit{i.e.} the map $(\varphi,\psi)\mapsto \langle \varphi,\opdb{\psi}{\mu} \rangle:=\int_{\R}\opdb{\psi}{\mu}\varphi(x)\dd x$  satisfies  $\langle \varphi,\opdb{\psi}{\mu}\rangle=\langle\psi,\opdb{\varphi}{\mu}\rangle$),  which is a property we use in our approach. Relaxing this assumption has potentially significant consequences.

When $\mu$ is nonsymmetric, a privileged direction is structurally introduced, and this affects the speed of propagation, making its estimation much more involved.  However, when $\mu$ is a bounded Borel measure earlier results of Yagisita \cite{Yagisita2009,Yagisita2009a} or \cite{Coville2007d,Coville2008a,Fang2015} show that all the above theorems remain true in such a context. Lebesgue's theory combined with Helly's theorem are sufficient to construct travelling wave solutions in such a situation. When the measure $\mu$ is unbounded, the nature of the singularity at $0$ plays a much more important role in the definition of the operator itself, making the generalisation of the above theorems more delicate. Indeed,
when $\mu$ is strongly singular in the sense that $\ds{\int_{\R}|z|\mathds{1}_{\{|z|\le 1\}}(z)d\mu(z)=+\infty}$, $\D_\mu$ needs to be replaced by the operator 
$$\opdb{u}{\mu}(x):=\int_{\R}[u(x+z)-u(x)-z u'(x)\mathds{1}_{\{|z|\le 1\}}(z)
] \dd \mu(z).$$
Our approach is mainly based on finding \textit{a priori} estimates on the speed which in turn provide $\mathsf{H}^2_{\textrm{loc}}(\R)$ regularity estimates that are  used to pass to the limit in a weak formulation of the equation. From a thorough inspection of our proofs, one could check that our arguments can be adapted when $\mu$ has a finite first moment, i.e. $\int_{\R}|z| d\mu(z)<+\infty$.  In such a case, the speed and the asymmetry of the kernel are of the same order and a proper weak formulation can be written for which estimates in  $\mathsf{BV}(\R)$ are sufficient to ensure the convergence. Such a  $\mathsf{BV}$ type estimate can naturally be obtained for bounded monotone functions. For a generic nonsymmetric unbounded measure, the existence of travelling wave solutions is largely an open question. Nevertheless, we believe that our approach can be used  when more structure on the measure $\mu$ is assumed. For example, when the measure $\mu\sim \frac{1}{|z|^{1+2s}}$ when $|z|\to 0$ with $s>\frac{1}{2}$ or when the support of $\mu$ is bounded. In both situations, we can rely on some regularity theory \cite{Caffarelli2009,Farkas2001} that provides the necessary estimates  to pass to the limit in a properly defined weak formulation.

\medskip

The organisation of the rest of the paper is as follows. In \Cref{sec:estimates}, we derive several general estimates on solutions to \eqref{eq:TW} that are of particular interest and used in later sections. In \Cref{sec:keylem}, a crucial lemma is proved in the bistable framework. Given these two sections, we shall prove \Cref{thm:bi}, \Cref{thm:ign}, and \Cref{thm:mono} in \Cref{sec:bi}, \Cref{sec:ign} and \Cref{sec:mono}, respectively. We complement this paper with \Cref{sec:tool}, where we recall various comparison principles and related results that are used throughout our arguments, together with important properties of operators $\D_\mu$.

 \subsection*{Notations}
$\N,\Z,\R$ and $\mathbb{C}$  denote the usual sets of numbers and their $^{\star}$ counterparts are the corresponding spaces without $0$. Classical functional spaces of continuous or differentiable functions will be denoted by $\mathscr{C}(\R),\mathscr{C}^{k}(\R),\mathscr{C}^{\infty}(\R)$ whereas for $1\le p< \infty,\;\mathsf{L}^p(\R),$ stands for the classical Banach space of $p-$integrable functions with respect to the  Lebesgue measure. $\mathsf{L}^\infty(\R),$ denotes the space of bounded functions. Spaces of 
 H\"older continuous functions and Sobolev spaces will be respectively denoted  by $\mathscr{C}^{k,\alpha}(\R)$ with $\alpha\in(0,1)$ and by $\mathsf{W}^{k,p}(\R)$. The  notation $\mathsf{H}^{k}(\R)$ will be used for $\mathsf{W}^{k,2}(\R)$. The Schwartz space will be denoted by $\mathscr{S}(\R)$. Spaces of continuous functions tending to $0$ at both infinities and with compact support will be denoted by $\mathscr{C}_0(\R)$ and $\mathscr{C}_c(\R)$ respectively. 
  We denote by  $\mathscr{X}_{\mathrm{loc}}$ the local version of a functional space $\mathscr{X}$. For $f\in \mathsf{L}^{1}(\R)$, we denote by $\mathscr{F}[f]$ its Fourier transform, i.e. $\mathscr{F}[f](\xi):=\int_{\R}e^{-iz\xi}f(z) \dd z$. For any measure $\mu$,  $\text{Supp}(\mu)$ denotes  its  support, i.e $\text{Supp}(\mu):= \{z\in \R^\star\,|\, \forall r>0,\,   \mu(B_r(z))>0\}$. Finally, for any positive real numbers $a < b$ and any measure $\mu$, we denote $\mu_{a}:=\mathds{1}_{\{|z|\ge a\}} \mu$  and  $\mu_{a,b} := \mathds{1}_{\{a \leq |z|\le b \}} \mu$, and therefore by $\D_{\mu_a}$  respectively $\D_{\mu_{a,b}}[u]$ the operator $\D$ with respectively the measure $\mu_{a}$  and the measure $\mu_{a,b}$.

\section{Estimates and useful formulas}

\subsection{A priori estimates}\label{sec:estimates}

We state some \textit{a priori} estimates  on monotone solutions to \eqref{eq:TW}, that is,
\begin{align*}
&cu' +\opdb{u}{\mu} +f(u)=0\quad \text{ on } \; \R ,\\[5pt]
&\lim_{x\to -\infty} u(x)=1 \qquad \lim_{x\to +\infty} u(x)=0
\end{align*} 
when $\mu$ is a bounded measure.  
\begin{prop}\label{prop:esti}
Assume that $\mu$ is bounded, that is, $\mu(\R)<+\infty$. Let $(c,u) \in \R^\star \times \left( \mathscr{C}^1(\R)\cap\mathsf{L}^{\infty}(\R)\right)$ be a strong monotone solution to \eqref{eq:TW}. Then  $u\in \mathscr{C}^2(\R)$, and 
satisfies:
\begin{enumerate}
\item\label{1} $c \Vert u' \Vert_{\sfL^2(\R)}^2 = \int_{0}^1 f(s)\, \dd s$,
\item\label{2} $|c| \| u'\|_{\infty} \le 2 \mu(\R)+ \Vert f \Vert_{\mathsf{L}^{\infty}([0,1])}$, 
\item\label{3} $c^2 \| u''\|_{\infty} \le \left(2\mu(\R) + \Vert f\Vert_{\mathsf{W}^{1,\infty}([0,1])} \right)^2 $
\item\label{4} $\ds{|c|^3\|u''\|_{\sfL^2(\R)}^2  \leq \left( \Vert f' \Vert_{\mathsf{L}^{\infty}([0,1])} \right)^2 \left|\int_{0}^{1}f(s)\, \dd s\right|}$.
\item\label{5} If further $\mu$ has a first moment (i.e. $\ds{\int_{\R}|z|\,\dd \mu(z)<+\infty}$) then \\  $\ds{c=\int_{\R}f(u(x))\,\dd x =-\int_{\R}\int_{\R} (u(x+z)-u(x))^2\,\dd \mu(z)\dd x +2\int_{\R}f(u(x))u(x)\,\dd x}$.
\end{enumerate}
\end{prop}   

\begin{proof}[{\bf Proof of \Cref{prop:esti}}]
The regularity of $u$ is obtained by using the equation satisfied by $u$. Indeed, observe that since $u$ is bounded and $\mu$ is a bounded Borel measure, the function 
$$\mu\star u (\cdot) :=\int_{\R}u(\cdot+z)\,\dd \mu(z)$$ 
is well defined, bounded, and 
$$\opdb{u}{\mu} = \mu\star u -  \mu(\R) \, u .$$ 
In addition,  $\mu\star u (\cdot)\in \mathscr{C}^1(\R)$ since  $u\in \mathscr{C}^1(\R)$.  
Since $c \neq 0$ and $f(u),u$ and $\mu\star u (\cdot)$ are $ \mathscr{C}^1(\R)\cap \mathsf{L}^{\infty}(\R)$ functions  we deduce that  $u'\in \mathscr{C}^1(\R)$ and $u'\in\mathsf{L}^{\infty}(\R)$. Thus, since $u$ satisfies \eqref{eq:TW}, $u\in \mathscr{C}^2(\R)$.

For such a smooth monotone solution, we easily get \textit{a priori} estimates. Indeed, first observe that since $u$ has limits at $\pm \infty$, $u\in\mathsf{L}^{\infty}(\R)$ and monotone we have  $\lim_{x\to\pm \infty}u'(x)=0$  since $u'\in\sfL^1$ and $u'$ is uniformly continuous thanks to $u''\in \mathsf{L}^{\infty}(\R)$. Moreover $u'\in\mathsf{L}^2$.  Now by testing \eqref{eq:TW} on $u'$, we get
$$
-c \Vert u' \Vert_{\sfL^2(\R)}^2 = \int_{\R} u' (x)(\mu \star u)(x)\,\dd x + \frac{1}{2}\mu(\R) - \int_{0}^{1}f(s) \, \dd s,
$$
By integration by parts,
$$
\int_\R u'(x) (\mu \star u)(x)\,\dd x = [u(x) (\mu\star u)(x) ]^{+\infty}_{-\infty} - \int_{\R} u(x) (\mu \star u)'(x) = -\mu(\R) - \int_{\R} u(x) (\mu \star u')(x)\,\dd x.
$$
By using Fubini's theorem and performing the change of variables $y=x+z$ and $z\to-z$, we have, thanks to the symmetry of $\mu$, 
$$
 \int_{\R} u(x) (\mu \star u')(x)\,\dd x=\int_\R\int_{\R}u(y+z)u'(y) \,\dd \mu(-z)\dd y= \int_\R u'(y) (\mu \star u)(y)\,\dd y.$$
Therefore, 
$$
2\int_\R u'(x) (\mu \star u)(x)\,\dd x = -\mu(\R),
$$
and $\eqref{1}$ is  proved.

Estimate $\eqref{2}$ comes directly from an evaluation of the equation \eqref{eq:TW} on a sequence that approximates a maximum point of the function $u'$. Estimate $\eqref{3}$ comes from the same argument with $u''$, after differentiating \eqref{eq:TW}. 
To obtain $\eqref{4}$, remembering that  $u\in \mathscr{C}^2(\R)$ and $f\in \mathscr{C}^1([0,1])$, we differentiate the equation satisfied by $u$,   
$$
-cu''= \opdb{u'}{\mu}+f'(u)u'.
$$ 
As for $u'$ we can show that $u''\in\mathsf{L}^2(\R)$ by using that $u'\in \mathsf{L}^{\infty}(\R)$, $c\neq 0$ and the above equation. Now, by  testing   on $u''$ the latter equation,
$$
-c \Vert u'' \Vert_{\sfL^2(\R)}^2 = \int_{-\infty}^{+\infty}\opdb{u'}{\mu}(x)u''(x)\,\dd x + \int_{-\infty}^{+\infty} f'(u(x))u'(x)u''(x)\,\dd x.
$$
By the Cauchy-Schwarz inequality, we then get 
\begin{equation*}
|c| \Vert u'' \Vert_{\sfL^2(\R)}^2 \le  \left|\int_{-\infty}^{+\infty} \opdb{u'}{\mu}(x)u''(x)\,\dd x\right| +  \Vert f' \Vert_{\mathsf{L}^{\infty}([0,1])} \Vert u' \Vert_{\sfL^2(\R)} \Vert u'' \Vert_{\sfL^2(\R)}.
\end{equation*}
We estimate the first term of the r.h.s as above to obtain

 \begin{align*}
 |c|\|u''\|_{\sfL^2(\R)}^2 &\le   \Vert f' \Vert_{[0,1]}  \Vert u' \Vert_{\sfL^2(\R)} \Vert u'' \Vert_{\sfL^2(\R)},
 \end{align*}
Combining with \eqref{1}, we have
\begin{align*}
|c|\|u''\|_{\sfL^2(\R)}  \leq  \Vert f' \Vert_{[0,1]}  \sqrt{\frac{\left|\int_{0}^{1}f(s)\, \dd s\right|}{|c|}},
\end{align*}
from which the desired estimate follows.
Last, to obtain \eqref{5}, let us observe that by successively using a Taylor expansion of $u$ and Fubini's theorem, we have
$$\int_{-R}^R \opdb{u}{\mu}(x)\,\dd x=\int_{0}^1\int_{\R}z[u(R+tz)-u(-R+tz)]\dd \mu(z)\dd t.$$
Since $\mu$ has a first moment, the Lebesgue convergence theorem then implies that 
$$\lim_{R\to +\infty}\int_{-R}^R \opdb{u}{\mu}(x)\,\dd x= -\int_{\R}z\dd \mu(z)=0, $$
thanks to the symmetry of $\mu$. Thus, $c=\int_{\R}f(u(x))\,\dd x$. 
To obtain the second equality, we just have to observe that by multiplying  \eqref{eq:TW} by $u$ and integrating over $(-R,R)$, we get
$$\frac{-c}{2}\left[u^2(R)-u^2(-R)\right] = \int_{-R}^R\int_{\R} u(x)(u(x+z)-u(x))\dd \mu(z)\dd x +\int_{-R}^Rf(u(x))u(x)\,\dd x.$$
To complete the equality, we need to rewrite the first term of the r.h.s adequately. To do so let us observe that 
\begin{multline*}
 \int_{-R}^R\int_{\R} u(x)(u(x+z)-u(x))\dd \mu(z)\dd x=\int_{-R}^R\int_{\R} (u(x)-u(x+z))(u(x+z)-u(x))\dd \mu(z)\dd x \\ +\int_{-R}^R\int_{\R}\int_{0}^1 z u(x+z)u'(x+tz)\dd t \dd \mu(z)\dd x. 
\end{multline*}
By using the Fubini theorem and  integrating by parts the last integral with respect to $x$, we get 
 \begin{multline*}
 \int_{-R}^R\int_{\R} u(x)(u(x+z)-u(x))\dd \mu(z)\dd x=-\int_{-R}^R\int_{\R} [u(x+z)-u(x)]^2\dd \mu(z)\dd x \\ +\int_{\R}\int_0^1z[u(R+z)u(R+tz)-u(-R+z)u(-R+tz)]\,\dd t\dd \mu(z) -\int_{-R}^R\int_{\R}\int_{0}^1 z u'(x+z)u(x+tz)\dd t \dd \mu(z)\dd x. 
\end{multline*}
By sending $R\to +\infty$, using the symmetry of $\mu$, we end up with 
 \begin{multline*}
 \int_{\R}\int_{\R} u(x)(u(x+z)-u(x))\dd \mu(z)\dd x=-\int_{\R}\int_{\R} [u(x+z)-u(x)]^2\dd \mu(z)\dd x \\  -\int_{\R}\int_{\R}\int_{0}^1 z u'(x+z)u(x+tz)\dd t \dd \mu(z)\dd x. 
\end{multline*}
By  changing to the variables $y=x+tz $ and $s=1-t$ in the last integral we get
\begin{align*}
\int_{\R}\int_{\R}\int_{0}^1 z u'(x+z)u(x+tz)\dd t \dd \mu(z)\dd x &= \int_{\R}\int_{\R}\int_{0}^1 z u'(y+sz)u(y)\dd s \dd \mu(z)\dd y\\
&=\int_{\R}\int_{\R}[u(y+z)-u(y)]u(y)\dd \mu(z)\dd y,
\end{align*}
and thus

$$\int_{\R}\int_{\R} u(x)(u(x+z)-u(x))\dd \mu(z)\dd x=-\frac{1}{2}\int_{\R}\int_{\R} [u(x+z)-u(x)]^2\dd \mu(z)\dd x.$$
\end{proof} 

\subsection{An integration by parts formula}

\begin{lemma}\label{lem:iden}
Assume that $\mu$ has a bounded support. Then for any $\delta>0$, $\varphi\in \mathscr{C}^{\infty}_c(\R)$ and $\psi \in \sfL_{\textrm{loc}}^{\infty}(\R)$, we have
$$
\int_{\R} \opdb{\psi}{\mu_\delta} (x)\varphi(x)\,\dd x = \int_{\R} \opdb{\varphi}{\mu_\delta} (x)\psi(x)\,\dd x.
$$
If we further assume that $\mu$ is a bounded Borel measure, then the above equality also holds for $\delta=0$.
\end{lemma}

\begin{proof}[{\bf Proof of \Cref{lem:iden}}]
For $\psi \in \mathsf{L}_{\textrm{loc}}^{\infty}(\R)$,  by definition of $\mu_\delta$, we have $\opdb{\psi}{\mu_\delta}\in \mathsf{L}^{\infty}_{\mathrm{loc}}(\R)$ since $\mu$ has a bounded support. Therefore for any $\varphi\in \mathscr{C}^{\infty}_c(\R),$ we have $\opdb{\psi}{\mu_\delta} \varphi \in \mathsf{L}^1(\R)$  and by Fubini's theorem we have 
\begin{align*}
\I:=\int_{\R} \opdb{\psi}{\mu_\delta} (x)\varphi(x)\,\dd x 
&= \int_{\R}\int_{|z|\ge \delta} [\psi(x+z) -\psi(x)]\varphi(x)\,\dd \mu(z)\dd x,\\
&=\int_{|z|\ge \delta} \left( \int_{\R}[\psi(x+z) -\psi(x)]\varphi(x)\,\dd x\right) \dd \mu(z).
\end{align*}
Let us now split  into two equal parts the integral in the right hand side, and changed the variables $y=x+z$ in the inner integral and $z\to -z$ in the outer integral in the last double integral, by using the symmetry of $\mu$ we then get  
\[
\I = \frac{1}{2} \int_{|z|\ge \delta} \int_{\R}[\psi(x+z) -\psi(x)]\varphi(x)\dd x \dd\mu(z) +\frac{1}{2} \int_{|z|\ge \delta} \int_{\R}[\psi(y) -\psi(y+z)]\varphi(y+z)\,\dd y \dd \mu(z),
\]
which by rearranging the terms rewrites and renaming the integration variable
\[\I=\frac12 \int_{|z|\ge \delta} \left(\int_{\R}[\varphi(x)-\varphi(x+z)]\psi(x+z)\,\dd x\right) \dd \mu(z)  + \frac12 \int_{|z|\ge \delta} \left(\int_{\R}[\varphi(x+z) -\varphi(x)]\psi(x)\,\dd x\right) \dd \mu(z). \]

By redoing the change of variables $y=x+z$ and $z\to-z$, one has,   
$$\int_{|z|\ge \delta} \left(\int_{\R}[\varphi(x)-\varphi(x+z)]\psi(x+z)\,\dd x\right) \dd \mu(z)= \int_{|z|\ge \delta} \left(\int_{\R}[\varphi(y+z)-\varphi(y)]\psi(y)\,\dd y\right) \dd \mu(z)$$
and thus, 
\begin{align*}
\int_{\R} \opdb{\psi}{\mu_\delta} (x)\varphi(x)\,\dd x = \int_{|z|\ge \delta} \left(\int_{\R}[\varphi(x+z) -\varphi(x)]\psi(x)\,\dd x\right) \dd \mu(z).
\end{align*}
We conclude by using again Fubini's theorem.
\end{proof}

\section{A key lemma}\label{sec:keylem}

In this section, we prove a key lemma for our analysis concerning solutions to \eqref{eq:TW}, when the measure has bounded support and the nonlinearity is a generic bistable function. In particular, we assume that $g$ is of class $\mathscr{C}^1(\R)$, with $g(0)=0, g'(0)<0$, and satisfies 
\begin{equation}\label{eq:bistable-g}
    \exists \,0<\theta<\gamma, \text{ such that } g(\theta)=g(\gamma)=0, \; g|_{(0,\theta)}<0,  \; g|_{(\theta,\gamma)}>0.
\end{equation}
Throughout this section, we denote by $\lambda$ a symmetric unbounded measure that satisfies
\begin{equation*}
    \int_{\R}\min(1,z^2)\,\dd \lambda(z) <+\infty, \qquad\text{and} \qquad \lambda(\{x\})=0 \quad \text{ for all }\quad x\in\R^{\star},
\end{equation*}
and with uniformly bounded support, i.e. $\text{Supp}(\lambda)\subset (-r_0,r_0)$ for some $r_0>0$. We associate to $\lambda$ a sequence of bounded Borel measures $(\lambda_n)_{n\in\N}$, defined by
\begin{equation*}
    \lambda_n:=\mathds{1}_{\R\setminus[-\eps_n,\eps_n]}\lambda, \text{ with } (\eps_n)_{n\in \N} \text{ such that } \lim_{n \to + \infty} \eps_n = 0.
\end{equation*}
Given all this, we have:
\begin{lemma}\label{prop:exist-bound-supported-measure}
Suppose that, for every $n\in\N$, $(c_n,m_{n}) \in \R_{+}^{\star} \times \mathscr{C}^{2}(\R)$ is a solution of 
\begin{align}\label{eq:bistable-m_n}  
\begin{dcases}
-c_{n} m_{n}' = \opdb{m_{n}}{\lambda_n} + g(m_{n}) \quad \textit{ everywhere }\\[5pt]
\lim_{x\to -\infty} m_{n}(x) = \gamma,\qquad \lim_{x\to +\infty} m_{n}(x) = 0.  
\end{dcases}
\end{align}
Assume further that $c_{n}$ tends to $0$ as $n$ tends to $+\infty$ and $m_n$ is monotone non increasing. Then, there exists a continuous function $\overline{m}$,  positive and monotone decreasing, such that, 
\begin{align*}
 &\int_{\R}\big[\overline{m}(x)\opdb{\varphi}{\lambda}(x) + \varphi(x)g(\overline{m}(x))\big] \dd x=0 \quad \textit{ for all } \varphi\in \mathscr{C}_{c}^{\infty}(\R) ,\\[5pt]
 &\lim_{x\to -\infty} \overline{m}(x) = \gamma,\qquad \lim_{x\to +\infty} \overline{m}(x) = 0.
\end{align*}
\end{lemma}

Since it is rather long, we decompose the proof into several steps. 
First, we obtain exponential decay estimates on $m_{n}$. From this exponential behaviour, we can define its first and second  anti-derivatives, $v_{n}$ and $w_{n}$, respectively. These are sequences of smooth functions that are uniformly bounded in $\mathscr{C}^{0,1}_{\mathrm{loc}}(\R)$ and  $\mathscr{C}^{1,1}_{\mathrm{loc}}(\R)$, respectively. This ensures that, up to extraction, we can  pass to the limit in an adequate sense \footnote{One checks that $\mathscr{C}^{1,1}_{\mathrm{loc}}(\R)\subset \text{Dom}(\D_\lambda)$ for any Lévy measure $\lambda$ with bounded support. This is due to the fact that for $\phi\in \mathscr{C}^{1,1}_{\mathrm{loc}}(\R)$, we have $|\phi(x+z)+\phi(x-z)-2\phi(x)| \lesssim z^2$ which is integrable against any L\'evy measure. This is therefore the natural function space in which to study strong convergence.}.

In a second step, we study the convergence of the sequences $(m_{n})_{n\in\N}, (v_{n})_{n\in\N}$ and $(w_{n})_{n\in \N}$ towards $\overline{m}$, $\overline{v}$ and $\overline{w}$, respectively and derive an equation satisfied by $\overline{w}$. Then, we argue that $\overline{m}$ has the expected limits at infinity. In a fourth step, we derive the equation satisfied by $\overline{v}$, and in a fifth we recover the equation satisfied by $\overline{m}$. In a sixth  and last step, we show that $\overline{m}$ is continuous. 

\medskip

Before going to the proof, we introduce the following sequence of linear forms on $\mathsf{L}^{\infty}_{\textrm{loc}}(\R)$.
\begin{equation*}
\forall\, n\in\N,\quad  V\left[\lambda_n \right]: \varphi \in {\mathsf{L}^{\infty}_{\textrm{loc}}(\R) \mapsto  \int_{\R} \varphi(z) \, z^2 \, \dd \lambda_n(z)}.
\end{equation*}
From the definition of $\lambda$ and $(\lambda_n)_{n\in \N}$, we have the weak - $\star$ convergence of $V\left[\lambda_n \right]$ to $V[\lambda]$.  

\medskip

\noindent\textbf{\# Step one: some uniform decay estimates on $m_{n}$.} 

\medskip
Since $g'(0)<0$ and $\lambda_n$ has a bounded support for all $n \in \N$, there exists a unique $\varsigma_n \in \R_{+}^{\star}$ such that 
$$ \int_{\R}(e^{-\varsigma_n z} -1)\,d \lambda_{n}(z)-c_{n}\varsigma_n +\frac{g'(0)}{2}=0.$$
As a consequence, from \eqref{eq:bistable-g} and \eqref{eq:bistable-m_n}, using \Cref{thm:cp},  one may show that  $m_n(\cdot)\lesssim  e^{-\varsigma_n \bullet}$. Such decay estimate enables to define 
\begin{equation*}\label{def:bi-veps}
v_{n}(\cdot):=\int_{\bullet}^{+\infty}m_{n}(z)\,\dd z.
\end{equation*} 
After integrating \eqref{eq:bistable-m_n},   $v_{n}$ satisfies for all $x\in\R$ 
\begin{equation}\label{eq:bi-veps}
-c_{n} v_{n}'(x)=\opdb{v_{n}}{\lambda_{n}}(x) +\int_{x}^{+\infty}g(m_{n}(y))\,\dd y.  
\end{equation}
The term $\int_{x}^{+\infty}g(m_n)$ is well defined since $\frac{g(m_n)}{m_{n}}$ is bounded and  $m_{n}$ decays exponentially fast.
Similarly, since one has  $ v_{n} \lesssim e^{-\varsigma_nz}$, the function  $$w_{n}(\cdot):=\int_{\bullet}^{+\infty}v_{n}(z)\,\dd z,$$ is thus also well defined. 
As well, after integrating \eqref{eq:bi-veps}, for all $x\in\R$, we have 
\begin{equation}\label{eq:bi-weps}
-c_{n} w_{n}'(x)=\opdb{w_{n}}{\lambda_{n}}(x)  +\int_{x}^{+\infty}\int_{y}^{+\infty}g(m_{n}(z))\,\dd z\dd y.
\end{equation}
To obtain a uniform decay estimate on $m_n$, $v_n$ and $w_n$ for large $n$, we only need to obtain a uniform estimate on $m_n$.

Define $p : z \mapsto \gamma e^{-\overline{\varsigma} \, z}$, for some $\overline{\varsigma}$ to be chosen later. The symmetry of $\lambda_n$ and the non negativity of $c_{n}$ lead to
\begin{align*}
  c_{n}p'+\opdb{p}{\lambda_n} +\frac{g(m_n)}{m_{n}}p  &=    \left(\int_{\R}(e^{-\overline{\varsigma} z} -1)\,\dd \lambda_{n}(z)-c_{n}\overline{\varsigma} +\frac{g(m_n)}{m_{n}}\right)p,\\
&\le \left(\int_{\R}[\cosh(\overline{\varsigma} z) -1]\,\dd \lambda_{n}(z) +g'(0)+ \left[\frac{g(m_n)}{m_{n}} -g'(0)\right]\right)p,\\
& \le \left(\int_{\R}[\cosh(\overline{\varsigma} z) -1]\,\dd \lambda(z) +g'(0) +\mathcal{I}_n +  \left[\frac{g(m_n)}{m_{n}} -g'(0)\right]\right)p,
  \end{align*}
 where $$ \mathcal{I}_n= \int_{\R}(\cosh(\overline{\varsigma} z) -1)\,\dd \lambda_{n}(z)- \int_{\R}(\cosh(\overline{\varsigma} z) -1)\,\dd \lambda(z).$$
On the one hand, since $V\left[\lambda_n \right]$ converges weak - $\star$ to  $V[\lambda]$ and the map $z \mapsto z^{-2} \left( \cosh(\overline{\varsigma} z)-1\right)$ is locally bounded, $|\mathcal{I}_n|\le \frac{\left| g'(0) \right|}{4}$ for $n$ large enough. On the other hand, since $g(0)=0$ and $g\in \mathscr{C}^1([0,\gamma])$, there exists $s_0\le \theta$ such that for all $s\in (0,s_0)$,
$$ \left|\frac{g(s)}{s} -g'(0)\right|\le \frac{\left\vert g'(0) \right\vert}{4} = -\frac{g'(0)}{4}.$$ 
Thus, by normalising $m_{n}$ by $m_{n}(0)=s_0$, since $m_{n}$ is monotone, we have for all $x\ge 0$, $m_{n}(x) \leq s_0$, and thus
$$ \left|\frac{g(m_n)}{m_{n}} -g'(0)\right|\le -\frac{g'(0)}{4},$$
on that set. 
Consequently,
\begin{align*}
  c_{n}p'+\opdb{p}{\lambda_{n}} +\frac{g(m_n)}{m_{n}}p & \le \left(\int_{\R}[\cosh(\overline{\varsigma} z) -1]\, \dd\lambda(z)+\frac{  g'(0) }{2} \right)p,
  \end{align*}
Observe that 
$$
h : \varsigma\mapsto \int_{\R}[\cosh(\varsigma z) -1]\,\dd \lambda(z) + \frac{g'(0)}{2},
$$ 
is well defined (since $\lambda$  has a bounded support), smooth, increasing and convex. Moreover, since  $h(0)= \frac{g'(0)}{2}<0$, there exists a unique $\overline{\varsigma} >0$ such that $h(\overline{\varsigma}) = 0$.

We conclude using \Cref{thm:cp}. Indeed, for $n \in \N$  large enough, on $\R_{+}$, $\frac{g(m_n)}{m_{n}}\le 0$ and 
\[ c_{n}p'+\opdb{p}{\lambda_{n}}+\frac{g(m_n)}{m_{n}}p  \le 0,
\]
together with $m_{n}\le p$ on $\R_{-}$. We get that $m_{n} (z) \le \gamma e^{-\overline{\varsigma} z}$ for all $z\in \R$. By definition of $v_n$ and $w_n$, they  have a similar asymptotic behaviour, that is
\begin{equation}\label{eq:exp-decay}
m_{n} (z) \le p:= \gamma e^{-\overline{\varsigma} z}\qquad \text{and}\qquad v_n \le \frac{1}{\overline{\varsigma}} p =\frac{\gamma}{\overline{\varsigma}}e^{-\overline{\varsigma} z} \qquad \text{and} \qquad w_n  \le \frac{1}{\overline{\varsigma}^2}p = \frac{\gamma}{\overline{\varsigma}^2} e^{-\overline{\varsigma} z}. 
\end{equation}

Equipped with these decay estimates and having the functions $v_{n},  w_{n}$ at hand, we are now in position to construct a continuous solution $\overline{m}$.

\medskip

\noindent\textbf{\# Step two: Convergence of the sequences $(m_{n})_{n\in\N}, (v_{n})_{n\in\N}$ and $(w_{n})_{n\in\N}$.} 

\medskip
Recall that we have normalised $m_{n}$ by $m_{n}(0)=s_0$ at the previous step of the proof.  By definition of $w_{n}$, $w_{n}''=m_{n}\ge 0$   is uniformly bounded since for all $n\in\N$, $0\le m_{n}<\gamma$.  Since for large $n \in \N$, $m_{n} \leq p$, $v_{n} \lesssim p$ and $w_{n} \lesssim p$, $w_{n}$ and $w_{n}'$  are locally uniformly bounded as well. 
By a diagonal extraction argument, $(w_{n})_{n \in\N }$ converges along subsequences to $\overline{w}$ in $\mathscr{C}^{1,\alpha}_{\mathrm{loc}}(\R)$ for all $\alpha\in(0,1)$. Moreover, $\overline{w}$ is non negative, non increasing and convex. In addition,  since $w_{n}''=m_{n}$ is monotone, up to extracting a subsequence, using Helly's theorem, $w_{n}''$ converges pointwise to some non negative  function $\overline{m}$. To enhance legibility, we omit extractions in the sequel. 

The dominated convergence theorem guarantees pointwise convergence of $(v_n)_{n\in\N}$ towards a function $\overline{v}$.
As well since $\overline{m}\in \mathsf{L}^{\infty}(\R)$ and $\int_{\R}z^2\dd \lambda(z)<+\infty$, 
the dominated convergence theorem guarantees the following pointwise convergence
\[\lim_{n\to+\infty}  \opdb{w_{n}}{\lambda_n}=\opdb{\overline{w}}{\lambda}.\]
Indeed,  by definition of $\lambda_n$, we have 
\[
    \opdb{w_{n}}{\lambda_n}(\cdot) =\int_{0}^1\int_{0}^1\int_{\R} s m_{n}(\cdot+s\nu z)\mathds{1}_{\R\setminus [-\eps_n,\eps_n]}(z)z^2\,\dd \nu \dd s\dd \lambda(z).
\]
Since $m_{n}\mathds{1}_{\R\setminus [-\eps_n,\eps_n]}$ converges to $\overline{m}$ pointwise, the convergence follows.

Moreover, note that since $m_{n} \leq p$, $v_{n} \lesssim p$,  $m_{n},v_{n}$ are uniformly bounded in $\mathsf{L}^{1}((r,+\infty))$  for all  $r\in \R$. This added to the fact that $s \mapsto\frac{g(s)}{s}$ is bounded and  $m_{n}$ converges to $\overline{m}$ pointwise, by Lebesgue's dominated convergence theorem, we get, 
$$
\lim_{n \to +\infty}\int_{\bullet}^{+\infty}\int_{y}^{+\infty}g(m_{n}(z))\,\dd z\dd y=\int_{\bullet}^{+\infty}\int_{y}^{+\infty}g(\overline{m}(z))\,\dd z\dd y,$$
pointwise.

Passing to the limit in \eqref{eq:bi-weps}, $\overline{w}$  satisfies for all $x\in\R$,
\begin{equation}\label{eq:bi-w}
\opdb{\overline{w}}{\lambda}(x) +\int_{x}^{+\infty}\int_{y}^{+\infty}g(\overline{m}(z))\,\dd z\dd y=0.  
\end{equation}

\medskip

\noindent\textbf{\# Step three: $\overline{m}$ has the expected limit at $\pm\infty$.}

\medskip

By the pointwise convergence of $\overline{m}$, we have  $\overline{m}(0)=s_0$ and $0\le \overline{m}\le \gamma$. If there exists $x_0$ such that $\overline{m}(x_0)=0$, by monotonicity of $\overline{m}$, $\overline{m}(x)=0$ for all $x\ge x_0$. We can thus define $x_1:=\inf\{x\ge 0, \overline{m}(x)=0\}$. By definition of $\overline{v}$ and $\overline{w}$, we get $v(x_1)=w(x_1)=0$ and thus at this point $\overline{w}$ satisfies  
$$
\int_{\R}[w(x_1+z)-w(x_1)]\,\dd \lambda(z)=0.
$$
Therefore, $w(x)=0$ for all $x\in x_1+(-r_0,r_0)\cap \text{Supp}
(\lambda)$ and in particular for some $x'<x_1$. Going back to the definition of $\overline{w}$, this implies that $v(x)=0$ for all $x\ge x'$ and thus $\overline{m}(x)=0$ for $x\ge x'$, contradicting the definition of $x_1$. The function $\overline{m}$ is thus positive. 
On the other hand, by \eqref{eq:exp-decay} we have $\overline{m}$, $v$ and $\overline{w}$ are bounded by a multiple of $p=\gamma e^{-\overline{\varsigma} z}$, which implies that  $\lim_{+\infty} \overline{m}= \lim_{+\infty} \overline{v}=\lim_{+\infty} \overline{w}=0$.
It remains to show that  
$\lim_{-\infty} \overline{m}=\gamma.$

Since $\overline{m}$ is monotone and bounded, $\lim_{-\infty} \overline{m}=l^-$ for some  $l^{-}\in (s_0,\gamma]$.   
We first show by contradiction that $l^{-}>\theta$.
Assume thus that $s_0\le l^{-}\le \theta$. For later use, let us
  denote by $G$ the  function 
$$G :=\int_{\bullet}^{+\infty}g(\overline{m}(z))\,\dd z.$$

Recall that by definition of $s_0$, for all $0\le s\le s_0$, we have $g(s)\le \frac{3}{4}g'(0)s$, and thus
$G(0)\le \frac{3}{4}g'(0)v(0)<0$. The early definition of $g$ gives that $g(\overline{m})\le 0$ on $\R$ and therefore $G$ is increasing and negative on $\R$. From this, we get, on $\R_{-}$,
$$
\int_{\bullet}^{+\infty}G(y)\,\dd y=\int_{\bullet}^{0}G(y)\,\dd y +\int_{0}^{+\infty}G(y)\,\dd y\le \int_{\bullet}^{0}G(y)\,\dd y,
$$
which gives, recalling $G(\cdot) \leq G(0)$ on $\R_{-}$,
\begin{equation}\label{eq:lim-w1}
\lim_{x\to -\infty} \int_{x}^{+\infty}G(y)\,\dd y=-\infty.
\end{equation}

On the other hand, since $\overline{w}\in \mathscr{C}^{1,1}(\R)$  and  $ \overline w''= \overline m$ almost everywhere,  by applying  Taylor's theorem, we get for almost every $x\in\R$ and almost every $z\in\R$,
$$
\overline{w}(x+z)-\overline{w}(x)=\overline{w}'(x)z+\int_{0}^1\int_{0}^{1}sz^2 \overline{m}(x+s\nu z)\,\dd \nu \dd s.
$$
Using the symmetry of $\lambda$ and since $\overline{m}$ is bounded and  nonincreasing, we get, for almost every $x\in\R$
$$
\opdb{\overline{w}}{\lambda}(x) =\int_{\R} [\overline{w}( x +z)-\overline{w}(x)]\,\dd \lambda(z)=\int_{0}^{1}\int_{0}^1\int_{\R}sz^2\overline{m}( x +s\nu z)\,\dd \nu \dd s\dd \lambda(z) \le \frac{l^{-}}{2} \int_{\R}z^2\,\dd \lambda(z).
$$
Returning to \eqref{eq:bi-w}, it  follows that, for almost every $x\in\R$,
\[
 \int_{x}^{+\infty} G(y)\,\dd y\ge - \frac{l^{-}}{2} \int_{\R}z^2\, \dd \lambda(z). 
 \]
  contradicting \eqref{eq:lim-w1}. Thus $l^{-}>\theta$.

Knowing this, we again argue by contradiction and assume that $l^-\in (\theta,\gamma)$. By definition of $g$, since $l^-\in (\theta,\gamma)$, by the monotonicity of $\overline{m}$, there exist  $x_0 \in \R$ and $\tau>0$ such that, for all $x\le x_0$, we have $g(\overline{m}(x))\ge \tau$. Therefore, for all $y\le x_0$, $G$ satisfies
$$
G(y)=\int_{y}^{+\infty}g(\overline{m}(z))\,\dd z =\int_{y}^{x_0}g(\overline{m}(z))\,\dd z +G(x_0)\ge \tau (x_0 -y) + G(x_0).
$$
For $x< x_0$, we obtain 
\begin{align*}
\int_{x}^{+\infty}G(y)\,\dd y&=\int_{x}^{x_0}G(y)\,\dd y + \int_{x_0}^{+\infty}G(y)\,\dd y,\\
&\ge \frac{\tau}{2}(x_0-x)^2+G(x_0)(x_0-x) + \int_{x_0}^{+\infty}G(y)\,\dd y,
\end{align*}
from which we get 
\begin{equation}\label{eq:lim-w2}
\lim_{x\to -\infty} \int_{x}^{+\infty}G(y)\,\dd y=+\infty.
\end{equation}
Recall that since $\overline{w}$ is convex,
$$
\opdb{\overline{w}}{\lambda}(\cdot) =\int_{\R}[\overline{w}(\cdot +z)-\overline{w}(\cdot)]\,\dd \lambda(z) =  \frac{1}{2}\int_{\R}[\overline{w}(\cdot +z)+\overline{w}(\cdot -z)-2\overline{w}(\cdot)]\,\dd \lambda(z)\ge 0. 
$$
Thus from \eqref{eq:bi-w}, it  follows that,
$$
\int_{\bullet}^{+\infty}G(y)\,\dd y = -\opdb{\overline{w}}{\lambda}(\cdot)\le 0
$$ 
which  contradicts \eqref{eq:lim-w2}. Hence, $l^-=\gamma$ and $\lim_{-\infty}\overline{m}=\gamma$.
  
\medskip

\noindent \textbf{\# Step four: deriving the equation satisfied by $\overline{v}$.}

\medskip

We now derive an equation for $\overline{v}$. On the one hand, since $m_{n}$ is of class $\mathscr{C}^2(\R)$ for all $n\in\N$, $v_{n}\in \mathscr{C}^{3}(\R)$ and since $m_n$ is uniformly bounded, $v_n$ is uniformly bounded in $\mathscr{C}^{0,1}_{\textrm{loc}}(\R)$. 

On the other hand, for any $\varphi\in \mathscr{C}_c^{\infty}(\R)$, by using \Cref{{lem:iden}} 
$$
 \int_{\R} \opdb{v_{n}}{\lambda_{n}}(x) \varphi(x)\,\dd x = \int_{\R} \opdb{\varphi}{\lambda_{n}}(x)v_{n}(x)\,\dd x. 
 $$
Since $\varphi\in \mathscr{C}^{\infty}_c(\R),$  $\opdb{\varphi}{\lambda}\in \mathscr{C}^{\infty}_{c}(\R)$ and by using that $V[\lambda_{n}]$ converges weak - $\star$ to $V[\lambda]$,
$$
\lim_{n \to +\infty}\opdb{\varphi}{\lambda_n} = \opdb{\varphi}{\lambda}.
$$

As a consequence, since $v_{n}$ converges to $\overline{v}$ in $\mathscr{C}^{0,\alpha}_{\textrm{loc}}(\R)$ with $\alpha\in (0,1)$, by Lebesgue's dominated convergence theorem,  
  $$
\lim_{n\to \infty}\int_{\R}\opdb{v_{n}}{\lambda_{n}}(x) \varphi(x)\,\dd x =\lim_{n\to \infty}\int_{\R} \opdb{\varphi}{\lambda_{n}}(x)v_{n}(x)\,\dd x=  \int_{\R} \opdb{\varphi}{\lambda} \overline{v}(x)\, \dd x. 
 $$
Moreover, we can check that for all $\varphi\in \mathscr{C}^{\infty}_c(\R),$
$$
\lim_{n \to +\infty }\int_{\R}\varphi(x)\int_{x}^{+\infty}g(m_{n}(z))\,\dd z\,\dd x = \int_{\R}\varphi(x)\int_{x}^{+\infty}g(\overline{m}(z))\,\dd z\,\dd x.
$$
Therefore,   $\overline{v}$ satisfies  for all $\varphi \in \mathscr{C}^{\infty}_c(\R)$,
$$
\int_{\R} \left(\opdb{\varphi}{\lambda}(x)\overline{v}(x) + \varphi(x)\int_{x}^{+\infty}g(\overline{m}(z))\,\dd z\right)\,\dd x =0. 
$$
We shall estimate more precisely the first part of the integral. Since for all $\delta>0$, 
$$
\int_{\R} [\varphi(\cdot+z) - \varphi(\cdot)] \, \dd \lambda(z) = \int_{|z|\le \delta}  [\varphi(\cdot+z) - \varphi(\cdot)]  \, \dd \lambda(z)\\ +\int_{|z|\ge \delta}  [\varphi(\cdot+z) - \varphi(\cdot)] \, \dd \lambda(z).
$$
By using \Cref{lem:iden}, since $\overline{v}\in \mathsf{L}^{\infty}_{\textrm{loc}}(\R)$, we get for all $\delta>0,$

$$
\int_{\R} \opdb{\varphi}{ \lambda}(x)\overline{v}(x)\, \dd x = \int_{\R}\int_{|z|\le \delta}  [\varphi(x+z) - \varphi(x)]\overline{v}(x)\, \dd \lambda(z)\dd x+\int_{\R}\opdb{\overline{v}}{\lambda_\delta}(x)\varphi(x)\,\dd x. 
$$
As a consequence for all $\delta>0$ and $\varphi\in \mathscr{C}^{\infty}_{c}(\R)$, we have 
$$
\int_{\R} \int_{|z|\le \delta}  [\varphi(x+z) - \varphi(x)]\overline{v}(x)\, \dd \lambda(z)\dd x +  \int_{\R}\left(\opdb{\overline{v}}{\lambda_\delta}(x) + G(x)\right)\varphi(x)\,\dd x =0. 
$$

We shall refine our estimate by observing that since $\varphi$ is smooth, by using Taylor's theorem, 
\begin{align*}
\left|\int_{\R} \int_{|z|\le \delta}  [\varphi(x+z) - \varphi(x)]\overline{v}(x)\, \dd \lambda(z)\dd x\right|&= \left|\int_{\R} \int_{|z|\le \delta} \iint_{[0,1]^2}  tz^2 \varphi''(x+tsz)\overline{v}(x)\, \dd t \dd s\dd \lambda(z)\dd x\right|\\
&\le \frac{1}{2}\left(\sup_{x\in (-r_0,r_0)+supp(\varphi)}\overline{v}(x)\right)\|\varphi''\|_{\infty}|\text{Supp}(\varphi)| \int_{|z|\le \delta} z^2 \,\dd \lambda(z),
\end{align*}
and so, we get  for all $\varphi \in \mathscr{C}_c^{\infty}(\R)$,
$$ \lim_{\delta\to 0} \int_{\R}\left(\opdb{\overline{v}}{\lambda_\delta}(x) + G(x)\right)\varphi(x)\,\dd x =0.$$

To obtain an equation for $\overline{v}$, it remains to show that we can interchange the limit with the integral. We ensure this by showing some integrability condition on the function $\varphi\opdiff{\overline{v}}$ where  $\opdiff{\overline{v}}(x,z):=\frac{1}{2}\left(\overline{v}(x+z)+\overline{v}(x-z)-2\overline{v}(x)\right)$.   
Namely, we  claim
\begin{lemma}\label{lem:v-unif-int}
  $\varphi \,\opdiff{\overline{v}}$  is  integrable in $\mathsf{L}^{1}(\R\times\R,\dd x\otimes \dd \lambda)$.
\end{lemma}

\begin{proof}[{\bf Proof of \Cref{lem:v-unif-int}}]
   Let us first observe that  by Tonelli's theorem we have 
$$
\int_{\R}\int_{\R}|\varphi(x)\opdiff{\overline{v}}(x,z)|\, \dd \lambda(z)\dd x= \int_{\R}\int_{\R}|\varphi(x)\opdiff{\overline{v}}(x,z)|\,\dd x \dd \lambda(z).
$$
We will show that the right hand side is bounded. To do so, let us observe that  since $\varphi\in \mathscr{C}_c^{\infty}(\R)$, there exists $a,b \in \R$ such that $supp(\varphi)\subset[a,b]$ and thus 
$|\varphi\opdiff{\overline{v}}|\le\|\varphi\|_{\infty} |\opdiff{\overline{v}}|\mathds{1}_{[a,b]}$ and  
$$\int_{\R}|\varphi(x)\opdiff{\overline{v}}(x,z)|\,\dd x\le \|\varphi\|_{\infty}\int_{a}^b |\opdiff{\overline{v}}(x,z)|\,\dd x.$$
Let us now estimate $\ds{\int_{a}^b|\opdiff{\overline{v}}|(x,z)\dd x}$.
 Recall that $\overline v$ is the pointwise limit of $v_n$ and since $v_n\in \mathscr{C}^3(\R)$ and convex decreasing, we have, by using Taylor's expansion, 
    \begin{align*}
    \int_{a}^{b}|v_n(x+z)+v_n(x-z)-2v_n(x)|\dd x&=\int_{a}^b\left |\int_{0}^1[v_{n}'(x+tz)-v'_{n}(x-tz)]z\,\dd t\right |\dd x,\\
    &=\int_{a}^b\left |\int_{0}^1\int_{-1}^1v_{n}''(x+stz)tz^2\, \dd s \dd t\right |\dd x,\\
    &= \int_{a}^b\int_{0}^1\int_{-1}^1-m_{n}'(x+tsz)tz^2\, \dd s\dd t\dd x.
    \end{align*}
    Since $-m'_n \ge 0$, by using Tonelli's theorem and integrating in $x$ we then get 
    \[  \int_{a}^{b}|v_n(x+z)+v_n(x-z)-2v_n(x)|\dd x= z^2\int_{0}^1\int_{-1}^1t[m_n(a+stz)-m_{n}(b+stz)]\dd s\dd t\le \gamma z^2,
    \]
    where we use that $0\le m_n\le \gamma$.
     Now, by passing to the limit as $n\to +\infty$, since $v_n\in \sfL^{\infty}_{\textrm{loc}}(\R)$, by Lebesgue's theorem, we deduce that 
     $$\int_{a}^{b}|\opdiff{\overline{v}}(x,z)|\dd x=\frac{1}{2}\int_{a}^b|\overline{v}(x+z)+\overline{v}(x-z)-2\overline{v}(x)|\,\dd x\le \frac{\gamma}{2}z^2.$$
Incidentally, it yields  
$$  \int_{\R}|\varphi(x)\opdiff{\overline{v}}(x,z)|\,\dd x\le \frac{\gamma\|\varphi\|_{\infty}}{2}z^2, $$
and thus 
$$ \iint_{\R^2}|\varphi(x)\opdiff{\overline{v}}(x,z)|\,\dd x\dd \lambda(z)\le \frac{\gamma\|\varphi\|_{\infty}}{2}\int_{\R}z^2 \dd \lambda(z)<+\infty.$$
\end{proof}

With this lemma at hand, we can derive the equation for $\overline{v}$. 
Indeed, we can observe that by symmetry of $\lambda$,
$$\lim_{\delta\to 0}\int_{\R}\varphi(x)\opdb{\overline{v}}{\lambda_\delta}(x)\,\dd x = \lim_{\delta\to 0} \int_{\R}\varphi(x)\left(\frac{1}{2}\int_{|z|>\delta}[\overline{v}(x+z)+\overline{v}(x-z)-2\overline{v}(x)]\,\dd \lambda(z)\right)\,\dd x$$
and so 
$$\lim_{\delta\to 0}\int_{\R}\varphi(x)\opdb{\overline{v}}{\lambda_\delta}(x)\,\dd x = \lim_{\delta \to 0} \iint_{\R^2}\varphi(x)\opdiff{\overline{v}}(x,z)\mathds{1}_{\{|z|\ge \delta\}}(z) \, \dd \lambda(z)\dd x. $$
Since $\varphi\opdiff{\overline{v}}\in \mathsf{L}^{1}(\R^2,\dd x\otimes \dd \lambda)$, by Lebesgue's theorem, we 
have
$$\lim_{\delta \to 0} \iint_{\R^2}\varphi(x)\opdiff{\overline{v}}(x,z)\mathds{1}_{\{|z|\ge \delta\}}(z) \, \dd \lambda(z)\dd x=\iint_{\R^2}\varphi(x)\opdiff{\overline{v}}(x,z)\, \dd \lambda(z)\dd x, $$
meaning that
$$\lim_{\delta\to0} \int_{\R} \opdb{\overline{v}}{\lambda_\delta}(x)\varphi(x)\,\dd x=\int_{\R} \opdb{\overline{v}}{\lambda}(x)\varphi(x)\,\dd x.$$
Therefore, for all $\varphi\in \mathscr{C}_c^{\infty}(\R)$, 
$$
\int_{\R}\left(\opdb{\overline{v}}{\lambda}(x) +G(x)\right)\varphi(x)\,\dd x =0,
$$
which in turn, implies, by the du Bois-Reymond lemma, that
\begin{equation}
    \label{eq:bar-v}
\opdb{\overline{v}}{\lambda} +G =0,
\end{equation}
almost everywhere. 

\medskip

\noindent\textbf{\# Step five: Deriving an equation for $\overline{m}$.}

\medskip

We derive an equation for $\overline{m}$, starting from the above equation. Recall that for all $n \in \N$,  $v_{n}$ is a convex function since we have $v_{n}''=-m_{n}'$ is positive. Therefore, since $v_{n}$ converges to $\overline{v}$ pointwise, $\overline{v}$ is also a convex function and by using the definition of $\D_{\lambda}$, we have
$$ 
\textrm{P.V.}\left(\int_{\R}  [\overline{v}(\cdot+z) - \overline{v}(\cdot) ]\, \dd \lambda(z)\right)=\frac{1}{2}\left(\int_{\R} [\overline{v}(\cdot +z)+\overline{v}(\cdot-z) - 2\overline{v}(\cdot)]\, \dd \lambda(z)\right)\ge 0.
$$
As a consequence, since  $\overline{v}$ is a continuous convex nonincreasing function, the function  $\overline{v}$ is  differentiable almost everywhere. 
Recall also that, 
$$ \lim_{n \to +\infty} v_{n}= \lim_{n \to +\infty}\int_{\bullet}^{\infty}m_{n}(z)\,\dd z  = \int_{\bullet}^{+\infty}\overline{m}(z)\,\dd z, $$
and so, since $\overline{m}$ is uniformly bounded, one has  $\overline v\in \mathscr{C}^{0,1}(\R)$. In addition, since $\overline{m}$ is monotone and bounded it has at most a  countable set of discontinuities, thus $\overline v' = -\overline{m}$ \textit{a.e.}.  

We are now in position to derive an equation satisfied by $\overline{m}$. We first differentiate equation \eqref{eq:bar-v} in the sense of distributions. 
Multiplying \eqref{eq:bar-v} by $\varphi'\in \mathscr{C}^{\infty}_c(\R)$ and integrating, we obtain
$$
\int_{\R}\varphi'(x)\left(\opdb{\overline v}{\lambda}(x) +G(x)\right)\,\dd x=0.
$$

Since $\overline{m}$ is continuous almost everywhere, $G$ is differentiable for almost every point $x\in\R$ with $G'(x)=-g(\overline{m}(x))$. Therefore, by using integration by parts, since $G\in \mathsf{L}^{\infty}_{\textrm{loc}}(\R)$ we have 
$$
\int_{\R}G(x)\varphi'(x)\,\dd x=-\int_{\R}G'(x)\varphi(x)\,\dd x= \int_{\R}g(\overline{m}(x))\varphi(x)\,\dd x \qquad \text{for all }\quad \varphi\in \mathscr{C}_c^{\infty}(\R).
$$
On the other hand,  we have 
$$ 
 \lim_{\delta\to 0}\int_{\R}\opdb{\overline v}{\lambda_\delta}(x)\varphi'(x)\,\dd x = \int_{\R}\left(\lim_{\delta\to 0} \opdb{\overline v}{\lambda_\delta}(x)\right)\varphi'(x)\,\dd x,
$$
which by using integration by parts 
$$
\int_{\R}\opdb{\overline v}{\lambda_\delta}(x)\varphi'(x)\,\dd x=\int_{\R}\opdb{\overline{m}}{\lambda_\delta}(x)\varphi(x)\,\dd x,
$$
then yields
$$ 
 \lim_{\delta\to 0}\int_{\R}\opdb{\overline{m}}{\lambda_\delta}(x)\varphi(x)\,\dd x = \int_{\R}\left(\lim_{\delta\to 0} \opdb{\overline v}{\lambda_\delta}(x)\right)\varphi'(x)\,\dd x. 
$$
We deduce that $\overline{m}$ satisfies for all $\varphi\in\mathscr{C}^{\infty}_c(\R)$, 
$$ 
 \lim_{\delta\to 0}\int_{\R}\opdb{\overline{m}}{\lambda_\delta}(x)\varphi(x)\,\dd x + \int_{\R}g(\overline{m}(x))\varphi(x)\,\dd x=0. 
$$
which by using again \Cref{lem:iden}  implies 
$$ \lim_{\delta\to 0}\int_{\R}\overline m(x)\opdb{\varphi}{\lambda_\delta}(x)\,\dd x + \int_{\R}g(\overline{m}(x))\varphi(x)\,\dd x=0. $$
Hence, for all $\varphi\in \mathscr{C}^{\infty}_c(\R)$, we get 
$$ \int_{\R}\overline m(x)\opdb{\varphi}{\lambda}(x)\,\dd x + \int_{\R}g(\overline{m}(x))\varphi(x)\,\dd x=0. $$

\medskip
\noindent\textbf{\# Step six: Continuity of $\overline{m}$.}

\medskip

Let us prove that $\overline{m}\in \mathscr{C}(\R)$. To do so, let us observe that since 
 $\overline{v}$ satisfies \eqref{eq:bar-v}, then for $\rho_\tau$ a mollifier, we have,
$$\rho_\tau\star \opdb{\overline v}{\lambda}+\rho_\tau\star G=0.$$
Since, $\overline{v}$ is convex, by  Tonelli's theorem, we have  
 $$ \rho_\tau\star \opdb{\overline{v}}{\lambda} =\opdb{\rho_\tau\star\overline{v}}{\lambda},$$
which gives $$\opdb{\rho_\tau\star\overline{v}}{\lambda}+\rho_\tau\star G=0.$$
Since $\rho_\tau\star \overline{v}\in \mathscr{C}^{\infty}(\R)$, we can differentiate the above equation to obtain 
$$  \opdb{\rho_\tau\star\overline {m}}{\lambda}+\rho_\tau\star g(\overline{m})=0.$$
Since  $g\in \mathscr{C}^1(\R)$ and $\overline{m}\in\mathsf{L}^{\infty}(\R)$, we get 
\begin{equation}\label{eq:molli-bar-m-bound}
\|\opdb{\rho_\tau\star \overline{m}}{\lambda}\|_{\sfL^\infty(\R)}\le \sup_{s\in [0,\gamma]}|g(s)|:=K.
\end{equation}

Let us define $\phi:=\rho_\tau\star \overline{m}$. Recall that $\overline{m}\in \sf{BV}(\R)$, so there exists a positive measure $\Pi$ such that $$\overline{m}(\bullet)=\int_{\bullet}^{+\infty}\dd\Pi(z).$$ By definition of $\phi$ and $\overline{m}$ we have 
$$ \phi(x)=\int_{\R}\rho_\tau(-z) \left(\int_{x+z}^{+\infty}\, \dd \Pi(z) \right)\,\dd z,$$
and so we get
$$
     \phi(x+h)+\phi(x-h)-2\phi(x)= \int_{\R}\rho_\tau(-z)\left( \int_{x-h+z}^{x+z}\, \dd \Pi(y) -\int^{x+h+z}_{x+z}\, \dd \Pi(y)\right) \dd z 
$$
which may be rewritten
$$
\phi(x+h)+\phi(x-h)-2\phi(x)= \int_{\R}\int_{\R}\rho_\tau(-z)\left( \mathds{1}_{(x-h+z, x+z)}(y)- \mathds{1}_{(x+z, x+z+h)}(y)\right)\, \dd \Pi(y)\, \dd z. 
$$
Since $\mathds{1}_{(x-h+z, x+z)}(y)=\mathds{1}_{(y-x,y-x+h)}(z)$ and  $\mathds{1}_{(x+z,x+h+z)}(y)=\mathds{1}_{(y-x-h,y-x)}(z)$, by using Fubini's theorem, we get
$$
\phi(x+h)+\phi(x-h)-2\phi(x)= \int_{\R}\left(\int_{\R}\rho_\tau(-z) [\mathds{1}_{(y-x,y-x+h)}(z)- \mathds{1}_{(y-x-h,y-x)}(z)]\, \dd z\right)\, \dd \Pi(y). 
$$
After defining  $$J_\tau(r,h):=\int_{\R}\rho_\tau(-z)\left[ \mathds{1}_{(r,r+h)}(z)- \mathds{1}_{(r-h,r)}(z)\right]\, \dd z,$$
we observe that,
$$ \phi(x+h)+\phi(x-h)-2\phi(x)= \int_{\R}J_\tau(y-x,h)\, \dd \Pi(y).$$
By integrating in $h$  over $(\delta,+\infty)$ the above expression, we get 
$$\int_{\delta}^{\infty} (\phi(x+h)+\phi(x-h)-2\phi(x))\,\dd\lambda(h)= \int_{\delta}^{+\infty}\left(\int_{\R}J_\tau(y-x,h)\, \dd \Pi(y)\right)\,\dd\lambda(h).$$
By using Fubini's theorem, we get 
\begin{equation}\label{eq:JJ}
\int_{\delta}^{+\infty}\left(\int_{\R}J_\tau(y-x,h)\, \dd \Pi(y)\right)\,\dd\lambda(h)=\int_{\R}\left(\int_{\delta}^{+\infty}J_\tau(y-x,h)\,\dd\lambda(h)\right)\,\dd \Pi(y).    
\end{equation} 
We have
\begin{lemma}\label{lem:J}
 $\exists K_\tau>0$, such that 
$$\sup_{r\in\R,\delta>0} \left|\int_{\delta}^{+\infty} J_\tau(r,h)\,\dd\lambda(h)\right|<K_\tau$$    
\end{lemma}
\begin{proof}[{\bf Proof of \Cref{lem:J}}]
Recall that
$$
J_\tau(r,h)=\int_{r}^{r+h}\rho_\tau(-z)\,\dd z - \int_{r-h}^{r}\rho_\tau(-z)\,\dd z.
$$
So $J_\tau(r,0)=0$ for all $r\in\R$. In addition we have 
\begin{align*}
&\partial_hJ_\tau(r,h)= \rho_\tau(-r-h) -\rho_\tau(h-r)=-h\int_{-1}^1\rho_\tau'(-r+th)\,\dd t\\
&\partial_{hh}J_\tau(r,h)= -\rho_\tau'(-r-h)-\rho_\tau'(h-r).
\end{align*}
So 
$\ds{\left|\partial_hJ_\tau(r,h)\right| \le 2h \|\rho'_\tau\|_{\infty} }$ and for all $r\in\R,$ $\partial_hJ_\tau(r,0)=0$ and $|\partial_{hh}J_\tau(r,h)|\le 2\|\rho'_\tau\|_{\infty}$.
As a consequence, by Taylor-Lagrange's theorem, we have for all $r \in \R$,  $ |J_\tau(r,h)|\le  h^2\|\rho'_\tau\|_{\infty}$ and 
thus 
$$ \left|\int_{\delta}^{\infty}J_{\tau}(r,h)\,\dd\lambda(h)\right|\le  \|\rho'_\tau\|_{\infty}\int_{\delta}^{\infty}h^2\,\dd\lambda(h).$$
\end{proof}

By taking the limit $\delta\to0$ in \eqref{eq:JJ}, on the left-hand side, we get  $\opdb{\phi}{\lambda}$ since $\phi\in \mathscr{C}^{\infty}(\R)$ and  on the right-hand side, the double integral converges since $\Pi$ is a finite Radon measure. 
Thus, we get 
\begin{equation}\label{eq:bar-m-fonda}
\opdb{\phi}{\lambda}(x)=\int_{\R}\left(\int_{0}^{+\infty}J_\tau(y-x,h)\,\dd\lambda(h)\right)\, \dd \Pi(y).
\end{equation}

Assume now, by contradiction, that $\overline{m}$ is not continuous. Since $\overline{m}$ is monotone decreasing and uniformly bounded, $\overline{m} \in \mathsf{BV}(\R)$. Since $\overline{m}$ is monotone, its discontinuities are jump discontinuities. Let us denote by $x_0$ a jump point of  $\overline{m}$, again by using that $\overline{m}\in \mathsf{BV}(\R)$, we can rewrite $\overline{m}$ as follows:  
$$\overline{m} =\int_{\bullet}^{+\infty}\dd\widetilde{m}(z)+\alpha_0\mathds{1}_{(-\infty,x_0)} =\int_{\bullet}^{+\infty}\, \dd \Pi(z)$$ 
where $\alpha_0>0$  and $\widetilde{m}$ is a positive Radon measure such that  $\widetilde{m}(\{x_0\})=0$ and  $\Pi=\widetilde{m} +\alpha_0\delta_{x_0}$. By \eqref{eq:bar-m-fonda}, we deduce 
\begin{equation}\label{eq:bar-m-fonda2}
\opdb{\phi}{\lambda}(x)=\int_{\R} \left(\int_{0}^{+\infty}J_\tau(y-x,h)\,\dd\lambda(h) \right)\dd\widetilde{m}(y)+\alpha_0\int_{0}^{+\infty}J_\tau(x_0-x,h)\,\dd\lambda(h). 
\end{equation}
Let us define 
$$M_\tau:= \sup_{r\in\R} \int_{0}^{+\infty}J_\tau(r,h)\,\dd\lambda(h),$$ 
then we have 
$$M_\tau\ge \int_{0}^{+\infty}J_\tau(0,h)\,\dd\lambda(h)\ge \int_{0}^{+\infty} \left(\int_{0}^{h}\rho_\tau(-z)\dd z -\int_{-h}^{0}\rho_\tau(-z)\,\dd z\right)\,\dd\lambda(h). $$ 
Now, let us  choose $\rho_\tau$ such that $\text{Supp}(\rho_\tau)\subset(-\tau,0)$ and symmetric around $-\frac{\tau}{2}$.
For this type of mollifier,  we deduce from  the above relation that
$$
M_\tau\ge\int_{0}^{+\infty} \int_{0}^{h}\rho_\tau(-z)\dd z\,\dd\lambda(h)\ge \lambda((\tau,+\infty))\to +\infty \quad\text{ as }\quad \tau\to 0.
$$
Note that by using the change of variable $z\to z'+\frac{\tau}{2}$,  that $\rho_\tau(z'-\frac{\tau}{2})=\rho_\tau(-z'-\frac{\tau}{2})$ and changing variables $z\to -z$, we have
\begin{align*}
J_\tau\left(r+\frac{\tau}{2},h\right) 
&= \int_{r}^{r+h}\rho_\tau(-z-\frac{\tau}{2})\dd z - \int_{r-h}^{r}\rho_\tau(-z-\frac{\tau}{2})\dd z\\
&=\int_{r}^{r+h}\rho_\tau(z-\frac{\tau}{2})\dd z - \int_{r-h}^{r}\rho_\tau(z-\frac{\tau}{2})\dd z \\
&=\int_{-r-h}^{-r}\rho_\tau(-z-\frac{\tau}{2})\dd z-  \int^{-r+h}_{-r}\rho_\tau(-z+\frac{\tau}{2})\dd z.
\end{align*}
Going back to the variable $z\to z'+\frac{\tau}{2}$, we get 
\begin{align*}
J_\tau\left(r+\frac{\tau}{2},h\right) &=\int_{-r+\frac{\tau}{2}-h}^{-r+\frac{\tau}{2}}\rho_\tau(-z)\dd z-  \int^{-r+\frac{\tau}{2}+h}_{-r+\frac{\tau}{2}}\rho_\tau(-z)\dd z=-J_\tau\left(-r+\frac{\tau}{2},h\right).
\end{align*}
As a consequence, for all $r\in\R$, we have 
\begin{equation}\label{eq:esti-below-Jtau}
\int_{0}^{\infty}J_\tau(r,h)\dd \lambda(h)\ge - M_{\tau}.
\end{equation}
By definition of $M_\tau$, there exists $r_\tau$ such that $\int_{0}^{+\infty}J_\tau(r_\tau,h)\,\dd\lambda(h)\ge \frac{M_\tau}{2}$. Necessarily, $\lim_{\tau \to 0}|r_\tau| =0$. Indeed, if not $|r_\tau|\ge r_0>0$ up to a sequence, and    
we have the following contradiction as $\tau \to 0$,
\begin{align*}
    &M_\tau\le  2\int_{-r_\tau}^{+\infty} \int_{0}^{r_\tau+h}\rho_\tau(-z)\dd z\,\dd\lambda(h)\le 2\lambda((r_0,+\infty)) &\quad \text{ if } \quad & r_\tau\le -r_0\\
    &0<M_\tau\le -2\int_{0}^{+\infty} \left( \int_{r_\tau-h}^{r_\tau}\rho_\tau(-z)\,\dd z\right)\,\dd\lambda(h) \le 0 &\quad \text{ if } \quad &r_\tau\ge r_0,\; \tau\le \frac{r_0}{2}.
\end{align*}
Take now a sequence $\tau_n\to 0$ such that $|r_{\tau_n}|\to 0$ and consider $x_n=x_0-r_{\tau_n}$. Evaluating \eqref{eq:bar-m-fonda2} at $x_n$, we get 
\begin{align*}
\opdb{\phi}{\lambda}(x_n)&= \int_{\R}\left(\int_{0}^{\infty}J_{\tau_n}(y-x_0+r_{\tau_n},h)\,\dd\lambda(h)\right)\,\dd\widetilde{m}(y) + \alpha_0\int_{0}^{+\infty}J_{\tau_n}(r_{\tau_n},h)\,\dd\lambda(h),\\
&\ge  \frac{\alpha_0}{2}M_{\tau_n}+ \int_{\R}\left(\int_{0}^{\infty}J_{\tau_n}(y-x_0+r_{\tau_n},h)\,\dd\lambda(h)\right)\,\dd\widetilde{m}(y). 
\end{align*}
Last, we obtain some estimate on the remaining integral. For $\delta>0,$ let $I_\delta:=(x_0-\delta,x_0+\delta)$,  then  
\begin{multline*}
\int_{\R}\left(\int_{0}^{\infty}J_{\tau_n}(y-x_0+r_{\tau_n},h)\,\dd\lambda(h)\right)\,\dd\widetilde{m}(y)=\int_{I_\delta}\left(\int_{0}^{\infty}J_{\tau_n}(y-x_0+r_{\tau_n},h)\,\dd\lambda(h)\right)\,\dd\widetilde{m}(y)\\ + \int_{\R\setminus I_\delta} \left(\int_{0}^{\infty}J_{\tau_n}(y-x_0+r_{\tau_n},h)\,\dd\lambda(h)\right)\,\dd\widetilde{m}(y) 
\end{multline*}
Recall that $\widetilde{m}$ is a Radon measure such that $\widetilde{m}(\{x_0\})=0$, so we can find $\delta_0>0$ small such that  $\widetilde{m}(I_{\delta_0})\le \frac{\alpha_0}{4}$. By using \eqref{eq:esti-below-Jtau} we have 
\[
\int_{\R}\left(\int_{0}^{\infty}J_{\tau_n}(y-x_0+r_{\tau_n},h)\,\dd\lambda(h)\right)\,\dd\widetilde{m}(y)\ge -\frac{\alpha_0}{4}M_{\tau_n}+ \int_{\R\setminus I_{\delta_0}} \left(\int_{0}^{\infty}J_{\tau_n}(y-x_0+r_{\tau_n},h)\,\dd\lambda(h)\right)\,\dd\widetilde{m}(y). 
\]

For all $y\in\R\setminus I_{\delta_0}$, since $r_{\tau_n} \to 0$, there exists $n_0$ such that $|y-x_0+r_{\tau_n}|\ge  \delta_0-|r_{\tau_n}|\ge \frac{\delta_0}{2}$, for all $n\ge n_0$. Increasing  $n_0$ if necessary, we can also assume that $\tau_n\le \frac{\delta_0}{4}$.
By using the definition of  $J_{\tau_n}(r,h)$,

If $y-x_0+r_{\tau_n}< 0$:
\begin{align*}
\int_{0}^{\infty}J_{\tau_n}(y-x_0+r_{\tau_n},h)\,\dd\lambda(h) &=\int_{-y+x_0-r_{\tau_n}}^{+\infty} \int_{0}^{y-x_0+r_{\tau_n}+h}\rho_{\tau_n}(-z)\dd z\,\dd\lambda(h)\\
&\le \lambda(x_0-r_{\tau_n}-y,+\infty))\\
&\le \lambda\left(\left(\frac{\delta_0}{2},+\infty\right)\right).
\end{align*}

If $y-x_0+r_{\tau_n}> 0$:
\begin{align*}\int_{0}^{\infty}J_{\tau_n}(y-x_0+r_{\tau_n},h)\,\dd\lambda(h) &=-\int_{y-x_0+r_{\tau_n}-\tau_n}^{+\infty} \int_{0}^{y-x_0+r_{\tau_n}}\rho_{\tau_n}(-z)\dd z\,\dd\lambda(h)\\
& \ge -\lambda((y-x_0+r_{\tau_n}-\tau_n,+\infty))\\
&\ge  -\lambda\left(\left(\frac{\delta_0}{4},+\infty\right)\right).
\end{align*}
Hence, we finally obtain for $n\ge n_0$,
$$
\opdb{\phi}{\lambda}(x_n)\ge  \frac{\alpha_0}{4}M_{\tau_n}- \lambda\left(\left(\frac{\delta_0}{4},+\infty\right)\right) \int_{\R}\dd\widetilde{m}(y). 
$$
From which we get a contradiction with \eqref{eq:molli-bar-m-bound}. 
Consequently, $\overline{m}$ has no jump discontinuities and thus is continuous.

\section{Bistable type nonlinearities - Proof of \texorpdfstring{\Cref{thm:bi}}{thm:bi}}\label{sec:bi}

In this section, we prove \Cref{thm:bi}. As previously stated in the introduction, in this situation, Yagisita's work \cite{Yagisita2009a} shows that a travelling wave solution to \eqref{eq:main} exists when the measure is bounded i.e. $\mu(\R)<+\infty$. To complete the proof of \Cref{thm:bi}, we  only need to consider the case of an unbounded measure, i.e. $\mu$ such that $\mu(\R)=+\infty$. As mentioned in the introduction, let us consider the approximated bounded Borel measure $\mu_\eps := \mathds{1}_{\R\setminus [-\eps,\eps]} \; \mu$ and the associated approximated problem:   
 \begin{align*}
 &-c_\eps u^{\prime}_{\eps}= \opdb{u_\eps}{\mu_\eps} + f(u_{\eps}) \quad \textit{ a.e.},\\
&\lim_{x\to -\infty} u_\eps(x) = 1,\qquad \lim_{x\to +\infty} u_\eps(x) = 0. 
    \end{align*}
Thanks to Yagisita's existence result \cite{Yagisita2009a},  there exists a unique speed $c_\eps$ and a bounded monotone decreasing profile $u_\eps$, solution to the latter, which is $\mathscr{C}^1(\R)$ when $c_\eps\neq 0$. 

The proof goes into several steps. We start by showing some elementary regularity properties on the family $(c_\eps,u_\eps)_{\eps>0}$. Then, we obtain some uniform upper and lower bounds on $c_\eps$. Finally, we construct the front and prove uniqueness of its speed. 

\subsection{\textit{A priori} regularity}

Let us obtain first some  elementary regularity properties on the family $(c_\eps,u_\eps)_{\eps>0}$. Namely, we prove,

\begin{lemma}\label{lem:esti2}
Assume that $\mu$ is unbounded. Then there exists $0 < \eps_0 \leq 1$ such that for all $\eps \le \eps_0$, $c_{\eps} > 0$ and $u_\eps\in \mathscr{C}^2(\R)$.
\end{lemma}

\begin{proof}[{\bf Proof of \Cref{lem:esti2}}]
Since $\mu$ is unbounded, we have $$\lim_{\eps \to 0}\int_{\R}\,\dd \mu_{\eps}(z) =+\infty.$$ 
Thus, there exists $\eps_1 >0$ such that for all $\eps \le \eps_1$, 
$$j_\eps:=\int_{\R}\,\dd \mu_\eps(z)> \sup_{s\in[0,1]} |f'(s)|.$$
Therefore, for any $\eps\le \eps_1$, the function $f_\eps:s\mapsto -j_\eps s +f(s) $ is  decreasing. 

Note that for $\eps\le \eps_1$, $u_\eps$ is a continuous function. Indeed, since $f$ satisfies \eqref{eq:bistable}, if $c_\eps\neq 0$ then by \Cref{prop:esti}, $c_\eps>0$ and $u_\eps\in \mathscr{C}^2(\R)$. Otherwise $c_{\eps}=0$ and by \eqref{eq:ueps} we can see that,  
$$
\int_{\R}u_\eps(\cdot +z)\,\dd \mu_\eps(z) \equiv -f_\eps(u_{\eps}).
$$
Since $\mu$ is a diffuse measure and $u_\eps$ a bounded monotone function, $x\mapsto\int_{\R}u_\eps(x+z)\,\dd \mu_\eps(z)$ is continuous. As a consequence, $f_\eps(u_\eps)$ is  continuous  as well and  by  monotonicity and  regularity of $f_\eps$, so is $u_{\eps}$.

On the other hand, by using that $f$ is of type \eqref{eq:bistable}, we can find $\nu>0$ and $f_\nu < f$ a smooth function of type \eqref{eq:bistable} such that, 
$$
f_\nu(0)=f_\nu(\theta)=f_\nu(1-\nu)=0, \quad f'_\nu(0) < f'(0),\quad \int_{0}^{1-\nu}f_\nu(s)\, \dd s>0.
$$
For such a $f_\nu$, thanks to Yagisita's result \cite{Yagisita2009a}, there exists $(c_{\eps,\nu}, u_{\eps,\nu})$ (with $u_{\eps,\nu}$ monotone decreasing) such that 
\begin{align*}
&c_{\eps,\nu}u'_{\eps,\nu}+\opdb{u_{\eps,\nu}}{\mu_\eps} +f_\nu(u_{\eps,\nu})=0, \qquad\textit{ a.e., } \\[5pt]
&\lim_{x\to -\infty}u_{\eps,\nu}(x)=1-\nu \qquad \lim_{x\to +\infty} u_{\eps,\nu}(x)=0.
\end{align*}
Moreover, by the same argument as for $u_{\eps}$, $u_{\eps,\nu}$ is at least continuous for $\eps$ small, say $\eps\le \eps_2$, and $c_{\eps,\nu}\ge0$.  

We shall show that $c_\eps$ is positive for all $\eps\le \eps_0:=\min\{1,\eps_1,\eps_2\}$. If not, then $u_\eps$ and $u_{\eps,\nu}$ are continuous functions that must satisfy 
\begin{align*}
&\opdb{u_{\eps}}{\mu_\eps} +f(u_{\eps})= 0 , \\[5pt]
&\opdb{u_{\eps,\nu}}{\mu_\eps} +f(u_{\eps,\nu})>0 , \\[5pt]
&\lim_{x\to -\infty}(u_{\eps} -u_{\eps,\nu})(x) = \nu>0, \\[5pt]
&\lim_{x\to +\infty}(u_{\eps} -u_{\eps,\nu})(x) = 0.
\end{align*}
By using \Cref{thm:nonlin-cp}, we  then obtain that  $u_{\eps}(\cdot+\tau)>u_{\eps,\nu}(\cdot)$ for all $\tau \in \R$, from which we get  the contradiction 
$$0<u_{\eps,\nu}(0)\le \lim_{\tau \, \rightarrow \, +\infty }u_{\eps}(\tau)=0.$$ 
Hence, $c_\eps> 0$ for all $\eps\le \eps_0$.
 \end{proof}

\subsection{Upper bound on the speed $c_\eps$}\label{subsec:uppboundceps}

We  prove the following upper bound for the speed $c_\eps$:

\begin{prop}\label{prop:esti-upper}
Assume that $\mu$ is unbounded and let $f$ be of type \eqref{eq:bistable}. Let $(c_\eps,u_\eps)$ be a solution to \eqref{eq:ueps} and $\eps_0$ be given by \Cref{lem:esti2}. Then there exists  $\overline{c}>0$ such that for all $\eps \le \eps_0$, $c_\eps\le \overline{c}$.
\end{prop}

\begin{proof}[{\bf Proof of \Cref{prop:esti-upper}}]
Fix $\eps<\eps_0$, where $\eps_0$ is given in \Cref{lem:esti2}. For such an $\eps$, $u_\eps$ is at least $\mathscr{C}^2(\R)$  and  the function $U_\eps : (t,x) \mapsto u_\eps(x-c_\eps t)$ satisfies the following Cauchy problem
\begin{equation}\label{eq:bi-CP-eps}
    \begin{dcases}
\partial_t U_\eps= \opdb{U_\eps}{\mu_\eps} +f(U_\eps),\\[5pt]
U_\eps(0,\cdot)=u_\eps.        
    \end{dcases}
\end{equation}

We aim to construct a function $w^+$ such that, for all $\eps <\eps_0$, $w^+\ge U_\eps$ and such that for some $\overline{c}>0$:
 $$
 \lim_{t\to +\infty} \sup_{x \, > \, ct} \;w^+(t,x)\le\theta  \quad\text{ for all }\quad c> \overline{c}. 
 $$

From this spreading property, we can directly infer $c_\eps\le \overline{c}$. Indeed, if not, there exists $x_0 \in \R$ such that $\theta < u_\eps(x_0)$, which leads to
$$\forall t \in \R_{+}, \qquad \theta <u_\eps(x_0)\le w^{+}(t,x_0+c_\eps t)\le \sup_{\{x>c_\eps t\}}w^+(t,x),$$
and thus to the following contradiction
$$\theta< \lim_{t\to +\infty} \sup_{x \, > \, c_\eps t}w^+(t,x)\le \theta.$$

Such a function $w^+$ can be obtained using Chen's construction \cite{Chen1997}. Let $\zeta\in \mathscr{C}^{\infty}(\R,[0,1])$ be a fixed function having the following properties:
$$
\zeta(s)= \begin{dcases} 0 &\text{ if } s \le 0,\\
                            1 &\text{ if } s \ge 4,
             \end{dcases}
\qquad 0< \zeta'(s) < 1, \quad\text{for}\quad s\in(0,4),  \text{ and }\quad |\zeta''(s)|\le 1$$
Let us fix $\rho \in \left(0, \inf\left\{
\frac{\theta}{4}, \frac{1-\theta}{4}\right\}\right),$ and denote by $m_\rho$ and $\zeta'_{\textrm{min}}$  the positive constants
$$m_\rho:=\min_{s \in [\rho, \theta -\frac{\rho}{2}]} -f(s), \qquad \zeta'_{\textrm{min}}:= \min_{s\in \left[\zeta^{-1}(\rho), \zeta^{-1}(1-\frac{\rho}{2})\right]} \zeta'(s).$$

\noindent By definition of $\mu$ and $\mu_\eps$, we can choose $R_\rho$ such that for all $\eps\le \eps_0$
$$
\int_{|z|\ge R_\rho}\,\dd \mu_\eps(z)=\int_{|z|\ge R_\rho}\,\dd \mu(z)\le \frac{m_\rho}{8},
$$
 and let us denote by  $M_\rho$ the following quantity
 $$
 M_\rho:= \frac{1}{2}\int_{|z|\le R_\rho} z^2\dd \mu(z).
 $$
 Define also, the following constant 
$$
\sigma_0:= \frac{\sqrt{\theta^2+2M_\rho m_\rho} -\theta}{2M_\rho},\qquad
\overline{c}:= \frac{1+\|f\|_{\infty}+m_\rho}{ (1-\theta)\sigma_0 \zeta'_{min}}.
$$
Observe that $\sigma_0$ satisfies
\begin{equation}\label{eq:bi-supsol0}
M_\rho \sigma_0^2 + \theta\sigma_0 = \frac{m_\rho}{2}.
\end{equation}

Let us mimic Chen's construction \cite{Chen1997}. Consider the function  
\begin{equation}
w^+ : (t,x) \mapsto (1 + \rho)- [1-(\theta - 2\rho)e^{-\sigma_0 t} ]\zeta(\sigma_0(x - \bar ct)). \label{eq:bi-supsol-def}
\end{equation}
For any $t\ge 0$,  $w^+(t,\cdot)$ is  non increasing and takes its values in the interval $(\rho +(\theta -2\rho)e^{-\sigma_0 t}, 1+\rho)$. Observe that
\begin{align*}
&w^+\left(t,x_1(t):= \overline{c}t+\frac{1}{\sigma_0} \zeta^{-1}\left(\frac{\rho}{1-(\theta -2\rho)e^{-\sigma_0 t}}\right) \right)=1,\\
&w^+\left(t,x_\rho(t):= \overline{c}t+\frac{1}{\sigma_0} \zeta^{-1}\left(1-\frac{\rho}{2(1-(\theta -2\rho)e^{-\sigma_0 t})}\right) \right)=\frac{3\rho}{2}+(\theta-2\rho)e^{-\sigma_0 t}.
\end{align*}
We shall store below the values of $\partial_t w^+(t,x)$ and $\partial_{xx} w^{+}(t,x)$, for all $(t,x) \in \R_{+} \times \R$, for later use:    
\begin{align}
&\partial_t w^+(t,x)= \overline{c}\sigma_0  \zeta'(\sigma_0(x-\overline{c}t))[1- (\theta -2\rho)e^{-\sigma_0 t}]  -\sigma_0(\theta -2\rho)e^{-\sigma_0 t}\zeta(\sigma_0(x-\overline{c}t)) \label{eq:bi-supsol1}\\
&\partial_{xx}w^+(t,x) = -\sigma^2_0[1-(\theta - 2\rho)e^{-\sigma_0 t} ]\zeta''(\sigma_0(x-\bar c t)). \label{eq:bi-supsol2}
\end{align}
 
From the monotonic behaviour of $w^+(t,\cdot)$ and the definition of $x_1(t)$, we have that for all $\eps>0$, $t\ge 0$  and  $x\le x_1(t)$, 
$$  w^+(t,x)\ge 1>U_\eps(t,x).$$  
In addition, since $\lim_{x\to -\infty}w^+(0,x)=1+\rho$ and $\lim_{x\to +\infty}w^+(0,x)=\theta -\rho$, up to a spatial shift of $u_\eps$, we have $w^+(0,\cdot)>u_\eps(\cdot)=U_\eps(0,\cdot)$.

In order to obtain $w^+\ge U_\eps$ by means of the comparison principle recalled in \Cref{thm:pcp},  we only need to check that   $w^+$ is a supersolution to \eqref{eq:bi-CP-eps} for $t\ge 0$ and $x\ge x_1(t)$, given that $w^+(0,\cdot)\ge u_\eps(\cdot)$.

\begin{lemma}\label{lem:supersolwplus}
The function $w^+$ is a supersolution to \eqref{eq:bi-CP-eps} for $t\ge 0$ and $x\ge x_1(t)$.
\end{lemma}

\begin{proof}[{\bf Proof of \Cref{lem:supersolwplus}}]
We first estimate $\opdb{w^+}{\mu_\eps}(t,x)$ in this region. Since $w^+\in \mathscr{C}^2(\R)$, by using the symmetry of $\mu$ and Taylor's theorem, 
\begin{align*}
\opdb{w^+}{\mu_\eps}(t,x) &= \int_{|z|\le R_\rho} [w^+(t,x+z)-w^+(t,x)]\,\dd \mu_\eps(z) + \int_{|z|\ge R_\rho} [w^+(t,x+z)-w^+(t,x)]\,\dd \mu(z)\\ 
                 &\le \iint_{[0,1]^2}\int_{|z|\le R_\rho}sz^2\partial_{xx}w^+(t, x +s\nu  z)\dd \nu \dd s \dd \mu_\eps(z) + 2\|w^+\|_{\infty} \int_{|z|\ge R_\rho} \dd \mu(z).
\end{align*}
By using \eqref{eq:bi-supsol2}, $|\zeta''|\le 1$ and the definition of $w^+$, we get for all $t>0,x\in \R$, 
\begin{equation}\label{eq:bi-supsol3}
\opdb{w^+}{\mu_\eps}(t,x)\le \frac{\sigma_0^2}{2}\int_{|z|\le R_\rho} z^2\,\dd \mu_\eps(z) + 4 \int_{|z|\ge R_\rho}\dd \mu(z)\le \sigma_0^2 M_\rho +  \frac{m_\rho}{2}.
\end{equation}

Let us now check that $w^+$ is a supersolution of \eqref{eq:bi-CP-eps} for $t\ge 0$ and $x\ge x_1(t)$. We treat separately the two situations: $x_1(t)\le x\le x_\rho(t)$ and $x\ge x_\rho(t)$.

\medskip

\noindent \# The case $x\ge x_\rho(t)$. In this situation, from \eqref{eq:bi-supsol1} and the definition of $w^+$, \eqref{eq:bi-supsol-def}, and since $\zeta'\ge 0$, we get, at $(t,x)$,

$$
\partial_t w^+ - \opdb{w^+}{\mu_\eps}-f(w^+)\ge -\sigma_0 \theta e^{-\sigma_0 t}- \opdb{w^+}{\mu_\eps} -f(w^+).   
$$
In addition, by using again \eqref{eq:bi-supsol-def}, for $x\ge x_\rho(t)$, we get that for all $t\ge0$, $\theta -\frac{\rho}{2}\ge  w^+\ge \rho$ and thus, we have 
$$ -f(w^+) \ge \min_{s \in (\rho, \theta -\frac{\rho}{2})} -f(s)= m_\rho>0.$$
By using \eqref{eq:bi-supsol3},  we get
$$
\partial_t w^+ - \opdb{w^+}{\mu_\eps}-f(w^+)\ge \frac{m_\rho}{2} -\sigma_0 \theta - \sigma^2_0 M_\rho = 0,   
$$
recalling the expression of $\sigma_0$.   

\medskip
    
\noindent \# The case $x_1(t)\le x\le x_\rho(t)$.
In this situation, observe that from the definition of $x_1(t)$ and $x_\rho(t)$ we have 

$$
\zeta^{-1}\left(1-\frac{\rho}{2}\right)\ge \zeta^{-1}\left(1-\frac{\rho}{2(1-(\theta -2\rho)e^{-\sigma_0 t})}\right)\ge \sigma_0 (x-\overline{c}t)\ge \zeta^{-1}\left(\frac{\rho}{1-(\theta -2\rho)e^{-\sigma_0 t}}\right)\ge  \zeta^{-1}\left(\rho\right).
$$
From the definition of $\zeta$, we have, for all $x_1(t)\le x\le x_\rho(t)$,
$$
 \zeta'(\sigma_0 (x-\overline{c}t))\ge \min_{s\in (\zeta^{-1}(\rho),\zeta^{-1}(1-\frac{\rho}{2}))} \zeta'(s)= \zeta'_{\min}>0,
$$
and from the definition of $\overline{c}$, by using \eqref{eq:bi-supsol1}  we deduce that  
$$
\partial_tw^+\ge \overline{c}\sigma_0(1-\theta) \zeta'_{\textrm{min}} -\sigma_0 \theta\ge 1+\|f\|_{\infty}+m_\rho -\sigma_0\theta.
$$
As a consequence we have,
$$
\partial_t w^+ - \opdb{w^+}{\mu_\eps}-f(w^+)\ge 1+\|f\|_{\infty}-f(w^+)+ m_\rho -\sigma_0 \theta -  \opdb{w^+}{\mu_\eps},   
$$
which by using \eqref{eq:bi-supsol3} and that $\|f\|_{\infty} \geq f(w^+) $  implies  
$$
\partial_t w^+ - \opdb{w^+}{\mu_\eps}(t,x)-f(w^+)\ge 1+ \frac{m_\rho}{2} -\sigma_0 \theta - M_\rho \sigma_0^2.
$$
Again by using \eqref{eq:bi-supsol0}, we then achieve
$$
\partial_t w^+ - \opdb{w^+}{\mu_\eps}-f(w^+)\ge 1\ge0.
$$
\end{proof}
The lemma being proved, the proposition follows.
\end{proof} 

\subsection{A positive uniform lower bound for $c_\eps$}

Next, we obtain a uniform positive lower bound for $c_\eps$. Namely, we prove 

\begin{prop}\label{lem:esti3}
Assume that $\mu$ is unbounded and diffuse. 
 Let $\eps_0>0$ be defined in \Cref{lem:esti2}, then  there exists $\underline{c}>0$ such that 
 $$ \forall \, \eps\in (0,\eps_0), \qquad c_\eps>\underline{c}.$$
\end{prop}

\begin{proof}[{\bf Proof of \Cref{lem:esti3}}]
The proof follows two steps. We first obtain an estimate from below by the speed of a problem defined with a measure with a bounded support $\mu_{0,r}:= \mathds{1}_{\{0\le |z|\le r\}}(z) \, \mu$ and then obtain a uniform estimate on this lower bound.

\medskip

\noindent \# \textit{Step one}: a first lower bound. 

For $r>\eps$, define the following measure, with bounded support, $\mu_{\eps,r}:= \mathds{1}_{\{\eps\le |z|\le r\}}(z)\mu=\mathds{1}_{\{0\le |z|\le r\}}\mu_\eps$ and the function $$f_{r}(s):=f(s)-s\cdot\mu(\{|z|\ge r\}), \qquad s \in \R_{+}.$$ 
Fix $r_0$ such that $f_{r_0}$ is bistable with three distinct zeros in $[0,1)$ (the smallest being zero and the largest being denoted by, say, $\gamma_{r_0}< 1$) and  $\int_{0}^{\gamma_{r_0}}f_{r_0}(s)\, \dd s > 0$.
Since $f_r$ converges to $f$ locally uniformly in $\mathscr{C}^1(\R)$, for  $r\ge \sup\{r_0,2\eps_0\}$, $f_r$ is of bistable type, with three zeros and such that $\int_{0}^{\gamma_r}f_{r}(s)\, \dd s > 0$. 
By Yagisita's existence result \cite{Yagisita2009a}, for such $\eps$ and $r$ there exists $(c_{\eps,r},u_{\eps,r}) \in \R_+ \times \mathsf{L}^{\infty}(\R)$ such that,
\begin{align*}
 &-c_{\eps,r} u^{\prime}_{\eps,r}= \opdb{u_{\eps,r}}{\mu_{\eps,r}} + f_r(u_{\eps,r}) \quad \textit{ a.e.},\\[5pt]
 &\lim_{x\to -\infty} u_{\eps,r}(x) = \gamma_{r},\qquad \lim_{x\to +\infty} u_{\eps,r}(x) = 0, \label{ -bc}
  \end{align*}
Indeed, $\mu_{\eps,r}$ is bounded and $f_{r}$ is a bistable function. Thanks to \Cref{lem:esti2}, we have $c_{\eps,r} \in \R_{+}^{\star}$ and $u_{\eps,r}\in \mathscr{C}^2(\R)$. 

Assume by contradiction that $c_\eps < c_{\eps,r}$. Since $f \geq f_r$ and $(c_\eps,u_\eps)$ satisfies \eqref{eq:ueps} we have
\begin{align*}
0&= c_\eps u^{\prime}_{\eps} + \int_{\R} [u_\eps(\cdot+z) - u_\eps(\cdot)] \, \dd \mu_\eps(z) +f(u_\eps)  \\
&= c_\eps u^{\prime}_{\eps} +f(u_\eps)+  \int_{[-r,r]} [u_\eps(\cdot+z) - u_\eps ] \, \dd \mu_\eps(z)+  \int_{\{|z|>r\}} [u_\eps(\cdot+z) - u_\eps] \, \dd \mu(z)\\ 
&\geq c_{\eps,r} u^{\prime}_{\eps} +f(u_\eps) + \opdb{u_\eps}{\mu_{\eps,r}}- u_\eps \, \mu( \lbrace |z| \geq r \rbrace)+\mu_{r}\star u_{\eps}\\
&\geq c_{\eps,r} u^{\prime}_{\eps} + \opdb{u_\eps}{\mu_{\eps,r}} + f_r(u_\eps).
\end{align*}

Therefore, $u_\eps$ and $u_{\eps,r}$ satisfy,
\begin{align*}
&-c_{\eps,r} u^{\prime}_{\eps}\ge \opdb{u_\eps}{\mu_{\eps,r}}  +f_r(u_\eps ), \\[5pt]
&-c_{\eps,r} u^{\prime}_{r,\eps}= \opdb{u_{\eps,r}}{\mu_{\eps,r}}   +f_r(u_{\eps,r} ),\\[5pt]
&\lim_{x\to -\infty} (u_{\eps}-u_{\eps,r})(x)=1-\gamma_{r}>0,\\[5pt]
&\lim_{x\to +\infty} (u_{\eps}-u_{\eps,r})(x)=0.
\end{align*}
By using \Cref{thm:nonlin-cp},  we deduce that 
$ u_{\eps}(\cdot + \tau)\ge u_{\eps,r}$, for all $\tau\in \R$ which yields the following contradiction
$$ 0=\lim_{\tau\to+\infty} u_\eps(\tau)\ge u_{\eps,r}(0)>0.$$
Hence $c_\eps \geq c_{\eps,r}$.

\medskip

\noindent \# \textit{Step two}: Estimate of $c_{\eps,r}$. 

Set $r=r_1>r_0$, with $\gamma_{r_1}>\gamma_{r_0}$. Let us prove that
\begin{equation*}
\liminf_{\eps \to 0} c_{\eps,r_1}>0.
\end{equation*}
Assume by contradiction that there exists an extraction of $(c_{\eps,r_1})_{\eps>0}$, denoted by $(c_{\eps_n,r_1})_{n \in \N}$ that converges to zero. Recall that $u_{\eps,r_1}$ satisfies
\begin{align}\label{eq:bistable-u2r0}
&-c_{\eps,r_1} u^{\prime}_{\eps,r_1}= \opdb{u_{\eps,r_1}}{\mu_{\eps,r_1}}   +f_{r_1}(u_{\eps,r_1} ).
\end{align}

On the one hand, observe that by definition,  $\mu_{0,r_1}$ and $\mu_{\eps_n,r_1}$ are respectively a Lévy measure with a bounded support and a sequence of approximated bounded Borel measures to $\mu_{0,r_1}$. Since, $f_{r_1}\in \mathscr{C}^1(\R)$ is a generic bistable function and  for all $n\in\N$, $u_{\eps_n,r_1}\in \mathscr{C}^2(\R)$, $c_{\eps_n,r_1}>0$ and $c_{\eps_n,r_1}$ tends to $0$, by \Cref{prop:exist-bound-supported-measure}, there exists a continuous function $\overline{u}$ such that $u_{\eps_n,r_1}$ tends to $\bar u$ pointwise and that satisfies
\[
    \lim_{x\to -\infty}\overline{u}(x)=\gamma_{r_1}\qquad \text{ and }\qquad  \lim_{x\to +\infty}\overline{u}(x)=0.
\]

On the other hand, from Yagisita's results \cite{Yagisita2009a}, the following problem:
\begin{align}\label{eq:bistable-s_eps-r_0-def}
&- s_{\eps} m^{\prime}_{\eps}= \opdb{m_{\eps}}{\mu_{\eps,r_1}}    + f_{r_0}(m_{\eps} ),\\
&\lim_{x\to-\infty}m_{\eps}(x)=\gamma_{r_0}, \qquad \lim_{x\to +\infty}m_{\eps}(x)=0\notag.
\end{align}
has a solution $(m_\eps,s_\eps)$ for $\eps$ small enough. Arguing as in \Cref{lem:esti2}, for all $\eps\le \eps_0$, $s_\eps> 0$ and $m_\eps\in \mathscr{C}^2(\R)$.
By using the definition of $f_{r}$, we have $f_{r_1} \geq f_{r_0}$ and thus by \Cref{thm:mono-c} we deduce that $s_\eps \le c_{\eps,r_1}$. As a consequence,  $$0\le \lim_{n \to +\infty} s_{\eps_n} \le \lim_{n \to +\infty} c_{\eps_n,r_1}=0.$$ 
As above, by \Cref{prop:exist-bound-supported-measure}, we get a continuous  $\overline{m}$ such that $m_\eps$ tends to $\overline{m} $ pointwise and such that 
\[
    \lim_{x\to -\infty}\overline{m}(x)=\gamma_{r_0}\qquad \text{ and }\qquad  \lim_{x\to +\infty}\overline{m}(x)=0.
\]
Let $\rho_\tau$ be a $\mathscr{C}^{\infty}_c(\R)$ mollifier, to be specified later on, then by convolving \eqref{eq:bistable-u2r0} and \eqref{eq:bistable-s_eps-r_0-def} and passing to the limit $\eps \to 0$, we then get  
\begin{align*}
&\opdb{\rho_\tau\star \overline{u}}{\mu_{0,r_1}} +\rho_\tau\star f_{r_1}(\overline{u})=0,\\
&\opdb{\rho_\tau\star \overline{m}}{\mu_{0,r_1}} +\rho_\tau\star f_{r_0}(\overline{m})=0,\\
&\lim_{x\to -\infty}\rho_\tau\star\overline{u}(x)=\gamma_{r_1}, \qquad \lim_{x\to +\infty} \rho_\tau\star\overline{u}(x)=0\\ 
&\lim_{x\to -\infty}\rho_\tau\star\overline{m}(x)=\gamma_{r_0}, \qquad \lim_{x\to +\infty} \rho_\tau\star\overline{m}(x)=0. 
\end{align*}
By definition of $f_r$ since $\gamma_{r_1}>\gamma_{r_0}$, we have $$\mu( \lbrace |z| \geq r_1 \rbrace)< a := \frac{\mu( \lbrace |z| \geq r_1 \rbrace)+\mu( \lbrace |z| \geq r_0 \rbrace)}{2}<\mu( \lbrace |z| \geq r_0 \rbrace)$$
and since $\overline{u}\in\mathscr{C}(\R)\cap \sfL^{\infty}(\R)$ and $f\in\mathscr{C}^1(\R)$ by \Cref{lem:tec-f}, we can find $\tau_1$ such that for all $\tau\le \tau_1$,
$$\rho_\tau\star f_{r_1}(\overline{u})\ge f(\rho_\tau\star \overline {u})- a \rho_\tau\star\overline{u}.$$
Similarly, since $\overline{m}\in \mathscr{C}(\R)\cap\sfL^{\infty}$, by \Cref{lem:tec-f},  we can find $\tau_2$ such that for all $\tau\le \tau_2$,
$$\rho_\tau\star f_{r_0}(\overline{m})\le f(\rho_\tau\star \overline{m})-a\rho_\tau\star \overline{m}.$$
For $\tau\le\min(\tau_1,\tau_2)$, we then obtain two continuous functions $\rho_\tau\star\overline{u}$ and $\rho_\tau\star\overline{m}$ that satisfy:  
\begin{align*}
&\opdb{\rho_\tau\star\overline{u}}{\mu_{0,r_1}}+f(\rho_\tau\star\overline{u})- a\rho_\tau\star\overline{u} \le 0, \\
&\opdb{\rho_\tau\star\overline{m}}{\mu_{0,r_1}}+f(\rho_\tau\star\overline{m})- a\rho_\tau\star\overline{m} \ge 0, \\
&\lim_{x\to -\infty}\rho_\tau\star\overline{u}(x)-\rho_\tau\star\overline{m}(x)=\gamma_{r_1} -\gamma_{r_0}> 0,\\
&\lim_{x\to +\infty} \rho_\tau\star\overline{u}(x)-\rho_\tau \star \overline{m}(x)=0.
\end{align*}
By  using \Cref{thm:nonlin-cp}, we then obtain  $\rho_\tau\star u(\cdot+h)\ge \rho_\tau\star\overline{m}$ for all $h \in\R$, which provides the following contradiction
$$0<\rho_\tau\star\overline{m}(0)\le \lim_{h\to+\infty}\rho_\tau\star \overline{u}(h)=0.$$ 
Hence, $\liminf_{\eps \to 0} c_{\eps,r_1}>0$ holds true and therefore there exists $\underline{c}$ such that for all $0<\eps \le \eps_0$,
$$
c_\eps\ge c_{\eps,r_1}>\underline{c}.
$$
\end{proof}

\subsection{Construction of the front}\label{ss:constr-front}
Equipped with \Cref{lem:esti2}, \Cref{prop:esti-upper} and \Cref{lem:esti3}, we shall prove \Cref{thm:bi}. To improve the presentation, we split our steps into several lemmas.

\begin{lemma}\label{lem:convbi}
    Up to a subsequence,  $(c_{\eps},u_{\eps})$ converges to $(c,u) \in \R_{+}^{\star} \times \mathscr{C}^{1,\alpha}_{\textrm{loc}}(\R)$, for any $\alpha \in \left(0,\frac{1}{2}\right)$. In addition,  $(c,u)$ satisfies, 
\begin{equation}\label{eq:bistable-u-pp}
    \opdb{u}{\mu} +  cu' + f(u) = 0\quad \textit{ almost everywhere}.
\end{equation}    
\end{lemma}

\begin{proof}[{\bf Proof of \Cref{lem:convbi}}]

Since $(c_\eps)_{\eps < \eps_0}$ is uniformly bounded and bounded away from zero, it converges, up to a subsequence, to $c>0$. Moreover, as a consequence of \Cref{prop:esti}, the sequence $(u'_{\eps})_{\eps < \eps_0}$ is uniformly bounded in $\mathsf{H}^1(\R)$. By setting $u_\eps(0)=\frac{\theta}{2}$, from a standard diagonal extraction procedure and the Sobolev embedding, up to subsequences,  $(u_{\eps})_{\eps < \eps_0}$  converges  in $\mathscr{C}^{1,\alpha}_{\textrm{loc}}(\R)$ for any $\alpha \in \left(0,\frac{1}{2}\right)$ to a function $u$ and $u'_\eps\to u'$ weakly in $\sf H^1$. In addition, $u(0)=\frac{\theta}{2}$ and $u$ is also monotone non increasing.

To show that $u$ satisfies \eqref{eq:TW}, we will proceed as in step four of the proof of \Cref{prop:exist-bound-supported-measure}. First, we pass to the limit in \eqref{eq:ueps} in the sense of distributions and show that $u$ is a weak solution to \eqref{eq:TW}. 
Then we show that $u$ satisfies the right boundary conditions.

Without any loss of generality we may assume that $\eps \le \eps_0\le1$.
Test \eqref{eq:ueps} on $\varphi \in \mathscr{C}_{c}^{\infty}(\R)$ to get,
$$
\int_{\R}\opdb{u_{\eps}}{\mu_\eps}(x)\varphi(x)\,\dd x + \int_{\R}\left( c_{\eps}u'_{\eps} + f(u_{\eps}(x))\right) \varphi (x)\,\dd x =0.
$$

In order to pass to the limit in the previous equation, we split $\mu_\eps$ into two parts:
\begin{equation*}
\mu_{\eps,1}:=\mathds{1}_{\{|z|<1\}}(z)\mu_\eps, \qquad \mu_{1}:=\mathds{1}_{\{|z|\ge 1\}}(z)\mu.
\end{equation*}
The above equality then decomposes,  
\begin{align*}
&\int_{\R} \opdb{u_{\eps}}{\mu_{\eps,1}}(x)\varphi(x)\,\dd x + \int_{\R}\left( \opdb{u_{\eps}}{\mu_{1}}(x) +  c_{\eps}u'_{\eps} + f(u_{\eps}(x))\right) \varphi (x)\,\dd x =0.
\end{align*}
Since $u_{\eps}$ converges to $u$ in $\mathscr{C}^{1,\alpha}_{\text{\textrm{loc}}}(\R)$ for all $\alpha \in (0,\tfrac12)$,  the second integral converges  and we have  
$$
\lim_{\eps \to 0} \int_{\R}\left( \opdb{u_{\eps}}{\mu_1}(x) +  c_{\eps}u'_{\eps} + f(u_{\eps}(x))\right) \varphi (x)\,\dd x= \int_{\R}\left( \opdb{u}{\mu_1}(x) +  cu' + f(u(x))\right) \varphi (x)\,\dd x.
$$  
To pass to the limit in the  first integral, that is $\ds{\int_{\R} \opdb{u_{\eps}}{\mu_{\eps,1}}(x)\varphi(x)\,\dd x}$, we will use arguments from {\bf $\#$ Step four} of the proof of \Cref{prop:exist-bound-supported-measure} above, combined with the technical lemma 
\begin{lemma}\label{lem:u-unif-int}
  For all $\varphi\in \mathscr{C}^{\infty}_{c}(\R)$, the function $\varphi \,\opdiff{u}$  is  integrable in $\mathsf{L}^{1}(\R\times\R,\dd x\otimes \dd \mu_{0,1})$.
\end{lemma}
Let us postpone for the moment the proof of \Cref{lem:u-unif-int}. By  using that $u_{\eps}$ and $u$ are $\mathscr{C}^{1,\alpha}_{\mathrm{loc}}(\R)$ functions, and reproducing the preceding arguments from {\bf $\#$ Step four} of the proof of \Cref{prop:exist-bound-supported-measure}, we can check that 
$$
\lim_{\eps\to 0}\int_{\R} \opdb{u_{\eps}}{\mu_{\eps,1}}(x)\varphi(x)\,\dd x= \lim_{\delta\to 0} \int_{\R} \opdb{u}{\mu_{\delta,1}}(x)\varphi(x)\dd x.
$$ 
Thanks to the symmetry of $\mu_{\delta,1}$, we also have 
$$\lim_{\delta\to 0}\int_{\R}\varphi(x)\opdb{u}{\mu_{\delta,1}}(x)\,\dd x = \lim_{\delta\to 0} \int_{\R}\varphi(x)\left(\frac{1}{2}\int_{\{\delta < |z|<1\}}[u(x+z)+u(x-z)-2u(x)]\,\dd \mu(z)\right)\,\dd x$$
and so 
$$\lim_{\delta\to 0}\int_{\R}\varphi(x)\opdb{u}{\mu_{\delta,1}}(x)\,\dd x = \lim_{\delta \to 0} \iint_{\R^2}\varphi(x)\opdiff{u}(x,z)\mathds{1}_{\{\delta< |z|< 1\}}(z) \, \dd \mu(z)\dd x. $$
Since by \Cref{lem:u-unif-int}, $\varphi\opdiff{u}\in \mathsf{L}^{1}(\R^2,\dd x\otimes \dd \mu_{0,1})$, by Lebesgue's theorem, we 
have
$$\lim_{\delta \to 0} \iint_{\R^2}\varphi(x)\opdiff{u}(x,z)\mathds{1}_{\{\delta < |z|< 1\}}(z) \, \dd \mu(z)\dd x=\iint_{\R^2}\varphi(x)\opdiff{u}(x,z)\mathds{1}_{\{|z|< 1\}}(z)\, \dd \mu(z)\dd x, $$
meaning that
$$\lim_{\delta\to0} \int_{\R} \opdb{u}{\mu_{\delta,1}}(x)\varphi(x)\,\dd x=\int_{\R} \opdb{u}{\mu_{0,1}}(x)\varphi(x)\,\dd x.$$
Thus, for all $\varphi\in \mathscr{C}^{\infty}_{c}(\R),$ we deduce that
$$
\int_{\R}\left( \opdb{u}{\mu}(x) +  cu' + f(u(x))\right) \varphi (x)\,\dd x =0
$$
which in turn, implies, by the du Bois - Reymond lemma, that
\begin{equation}
    \opdb{u}{\mu} +  c \, u' + f(u) = 0,
\end{equation}
almost everywhere. 

\begin{proof}[{\bf Proof of \Cref{lem:u-unif-int}}]
The proof follows the same line of argument used in \Cref{lem:v-unif-int}.  Let us first observe that  by Tonelli's theorem we have 
$$
\int_{\R}\int_{\R}|\varphi(x)\opdiff{u}(x,z)|\, \dd \mu_{0,1}(z)\dd x= \int_{\R}\int_{\R}|\varphi(x)\opdiff{u}(x,z)|\,\dd x \dd \mu_{0,1}(z).
$$
We will show that the right hand side is bounded. To do so, let us observe that  since $\varphi\in \mathscr{C}_c^{\infty}(\R)$, there exists $a,b \in \R$ such that $\text{Supp}(\varphi)\subset[a,b]$ and thus 
$|\varphi\opdiff{u}|\le\|\varphi\|_{\infty} |\opdiff{u}|\mathds{1}_{[a,b]}$ and therefore 
$$\int_{\R}|\varphi(x)\opdiff{u}(x,z)|\,\dd x\le \|\varphi\|_{\infty}\int_{a}^b |\opdiff{u}(x,z)|\,\dd x.$$
Let us now estimate $\ds{\int_{a}^b|\opdiff{u}(x,z)|\dd x}$.
 Recall that $u$ is the pointwise limit of $u_{\eps}$ and since $u_{\eps}\in \mathscr{C}^2(\R)$, we have by using Taylor's expansion and the Cauchy-Schwarz inequality 
    \begin{align*}
    \int_{a}^{b}|u_{\eps}(x+z)+u_{\eps}(x-z)-2u_\eps(x)|\dd x&=\int_{a}^b\left |\int_{0}^1[u_{\eps}'(x+tz)-u'_{\eps}(x-tz)]z\,\dd t\right |\dd x,\\
    &\le  \int_{0}^1\int_{-1}^1\left(\int_{a}^{b} |u''_\eps(x+tsz)| \dd x \right)tz^2\, \dd s\dd t\\
    &\le \sqrt{|a-b|}\|u''_\eps\|_{\mathsf{L}^2(\R)} z^2.
    \end{align*}
  By \Cref{prop:esti} and \Cref{lem:esti3}  for $\eps\le \eps_0$, we have 
  $$ \|u''_\eps\|_{\mathsf{L^2}}^2\le \frac{\sup_{s\in(0,1)}\|f'(s)\|^2\int_{0}^1f(s)\dd s}{\underline{c}^3}=:K^2.$$
By passing to the limit as $\eps\to 0$, since $u_\eps\in \sfL^{\infty}(\R)$,  by Lebesgue's theorem, we deduce that 
     $$\int_{a}^{b}|\opdiff{u}(x,z)|\dd x\le \frac{K\sqrt{|a-b|}}{2}\,z^2.$$
Incidentally, it yields  
$$  \int_{\R}|\varphi(x)\opdiff{u}(x,z)|\,\dd x\le\frac{K\|\varphi\|_{\infty}\sqrt{|a-b|}}{2}\,z^2, $$
and thus 
$$ \iint_{\R^2}|\varphi(x)\opdiff{u}(x,z)|\,\dd x\dd \mu_{0,1}(z)\le \frac{K\|\varphi\|_{\infty}\sqrt{|a-b|}}{2}\int_{\R}z^2 \dd \mu_{0,1}(z)<+\infty.$$

\end{proof}

\end{proof}

Let us now show that this equality \eqref{eq:bistable-u-pp} holds everywhere. 
    
\begin{lemma}\label{lem:bi-u}
    The weak solution  $(c,u)$ constructed in \Cref{lem:convbi} satisfies  
\begin{equation*}
    \opdb{u}{\mu} +  cu' + f(u) = 0\quad \textit{ everywhere}.
\end{equation*}    
\end{lemma}
\begin{proof}[{\bf Proof of \Cref{lem:bi-u}}]
    To do so, let us remark that since $u'\in \mathsf{H}^1(\R)$ and $f\in \mathscr{C}^{1}(\R)$, we have 
 $$\opdb{u}{\mu_{1}}+c u' +f(u) \in \mathscr{C}^{0,\alpha}_{\mathrm{loc}}(\R),$$ for $\alpha\in (0,\frac{1}{2})$.   
We only need to prove that the map $x\mapsto \opdb{u}{\mu_{0,1}}(x)$ is continuous. To prove the continuity of  $\opdb{u}{\mu_{0,1}}$, we use a standard localisation procedure combined with the Jacob-Schilling type characterisation  given in \Cref{lem:dom-mu}. 

By \Cref{lem:dom-mu}, we have $$\mathcal{R}_{\psi_{0,1}}:=\{v\in \mathscr{C}_0(\R)\cap \sfL^1(\R), \psi_{0,1}\mathscr{F}[v] \in \sfL^1(\R)\} \subset \{v\in \mathscr{C}(\R)\cap \text{Dom}(\D_{\mu_{0,1}})\,|\, \opdb{v}{\mu_{0,1}}\in \mathscr{C}(\R)\}$$ and  for $\varphi \in \mathcal{R}_{\psi_{0,1}}$,   
$$\opdb{\varphi}{\mu_{0,1}}= \mathscr{F}^{-1}[\psi_{0,1}\mathscr{F}[\varphi]],$$
where $\psi_{0,1}$ denotes the Lévy-Khintchine (Fourier type) symbol of $\D_{\mu_{0,1}}$.

To finish the proof, we need to find localisations of $u$  that belong to  $\mathcal{R}_{\psi_{0,1}}$. To do so, take $\chi\in \mathscr{C}_{c}^{\infty}(\R)$ a positive cut-off function that satisfies 
$$\chi(x):=
\begin{dcases}
1\quad \text{ for } |x|\le 2\\
0\quad \text{ for } |x|\ge 3,
\end{dcases}
$$
and define its translates $\chi_R:=\chi\left(\cdot-R\right)$. From this, we shall truncate $u$ and define $u_R:=u\chi_R$. By definition, $u_R\in \sf H^2(\R)$ and $u_R\equiv u$ in $[R-2,R+2]$. In addition,  since $\chi_R\in\mathscr{C}_c^{\infty}(\R)$ and $u\in \mathscr{C}^{1,\alpha}(\R)$, for some $C>0$, we have for all $x,z \in \R,\; |(\chi_R(x+z)-\chi_R(x))(u(x+z)-u(x))|\le Cz^2$ and  
\begin{align*}
\opdb{u_R}{\mu_{0,1}}
&= \Xi_R+  u\opdb{\chi_R}{\mu_{0,1}} + \chi_R \opdb{u}{\mu_{0,1}}, \label{eq:bi-reg-uR-main}
\end{align*}
with $$ \Xi_R := \int_{\R}[\chi_R(\cdot+z)-\chi_R(\cdot)][u(\cdot+z)-u(\cdot)]\,\dd \mu_{0,1}(z) .$$
Observe that since $u$   satisfies \eqref{eq:bistable-u-pp}, we have 
$$ \opdb{u}{\mu_{0,1}} = -\opdb{u}{\mu_{1}}-c u' -f(u)$$
almost everywhere and since $\chi_R\in \mathscr{C}_c^{\infty}(\R)$ and $cu'+\opdb{u}{\mu_{1}} +f(u) \in \sf H^1_{\mathrm{loc}}(\R)$,
we deduce that
$\chi_R\opdb{u}{\mu_{0,1}}\in \sf H^1(\R)$. 
Observe also that since $\chi_R\in \mathscr{C}^{\infty}_{c}(\R)$ then $\opdb{\chi_R}{\mu_{0,1}}\in \mathscr{C}_{c}^{\infty}(\R)$ and therefore $u\opdb{\chi_R}{\mu_{0,1}}\in \sf H^1(\R)$ since $u\in \mathscr{C}^{1,\alpha}(\R)$. Last, let us observe that  $\Xi_R$ also belongs to $\sf H^1(\R)$. Indeed, $\Xi_R$ may be rewritten as
$$
\Xi_R(\cdot)=\int_{\R}\int_{0}^{1}\int_{0}^1\chi'_R(\cdot+tz)u'(\cdot+sz)z^2\dd t\dd s \dd \mu_{0,1}(z) \in \mathscr{C}_{c}(\R),
$$
and so $\Xi_R\in \sfL^2(\R)$. Now, since $\chi'_R,u'\in \sf H^1(\R)$ and $z^2\mu_{0,1}$ is a finite Radon measure, by Lebesgue's theorem we can differentiate the above expression to obtain 
\begin{align*}
&\Xi'_R\\
&=\int_{\R}\int_{0}^{1}\int_{0}^1\left(\chi'_R(\cdot+tz)u'(\cdot+sz)\right)'z^2\dd t\dd s \dd \mu_{0,1}(z)\\
&=\int_{\R}\int_{0}^{1}\int_{0}^1\chi''_R(\cdot+tz)u'(\cdot+sz)z^2\dd t\dd s \dd \mu_{0,1}(z) + \int_{\R}\int_{0}^{1}\int_{0}^1\chi'_R(\cdot+tz)u''(\cdot+sz)z^2\dd t\dd s \dd \mu_{0,1}(z)
\end{align*}

Let us consider both integrals separately. 
Since $u'\in \sf H^{1}(\R)$, the first one belongs to $\mathscr{C}_{c}(\R)$ and so to $\sfL^2$. Whereas the last triple integral satisfies by using the Cauchy-Schwarz inequality with respect to the finite Radon measure $z^2\mu_{0,1}$ and Jensen's inequality
\begin{align*}
&\left(\int_{\R}\int_{0}^{1}\int_{0}^1\chi'_R(\cdot+tz)u''(\cdot+sz)z^2\dd t\dd s \dd \mu_{0,1}(z)\right)^2 \\
&\le K \int_{\R}\left[\int_{0}^1\int_{0}^1\chi_R'(x+tz)u''(x+sz)\dd t\dd s   \right]^2 z^2\dd \mu_{0,1}(z)\\
   &\le K \int_{\R}\int_{0}^1\int_{0}^1(\chi_R'(x+tz))^2(u''(x+sz))^2\dd t\dd s \,z^2\dd \mu_{0,1}(z)
\end{align*}
with $K:=\ds{\int_{\R}z^2\dd \mu_{0,1}(z)}$.
Since $\chi'_R\in \sf H^1(\R)$ and $u''\in \sfL^2(\R)$, the latter then implies that the last triple integral belongs also to $\sfL^2(\R)$ and thus $\Xi_R\in \sf H^1(\R)$.
Next, we prove that 
\begin{lemma}\label{lem:XI}
Almost everywhere, 
    \begin{equation*}
    \psi_{0,1} \mathscr{F}[u_R]= \mathscr{F}[\Xi_R] +\mathscr{F}[u\opdb{\chi_R}{\mu_{0,1}}] +\mathscr{F}[\chi_R\opdb{u}{\mu_{0,1}}]. 
\end{equation*}
\end{lemma}
 \begin{proof}[{\bf Proof of \Cref{lem:XI}}]
    For $0 < \delta <1$, a quick algebraic computation shows
    \begin{equation}
\opdb{u_R}{\mu_{\delta,1}}
= \Xi_{\delta,R}+  u\opdb{\chi_R}{\mu_{\delta,1}}+ \chi_R \opdb{u}{\mu_{\delta,1}}, \label{eq:bi-reg-uR-delta}
\end{equation}
with $$ \Xi_{\delta,R} := \int_{\R}[\chi_R(\cdot+z)-\chi_R(\cdot)][u(\cdot+z)-u(\cdot)]\,\dd \mu_{\delta,1}(z) .$$
As above, we can check that $\Xi_{\delta,R}, u\opdb{\chi_R}{\mu_{\delta,1}}$ and $\chi_R\opdb{u}{\mu_{\delta,1}}$ belong to $\sfL^2(\R)$ and since $u_R\in\sf H^2(\R)$ and the measure $\mu_{\delta,1}$ is a finite Borel measure, by taking the Fourier transform of \eqref{eq:bi-reg-uR-delta}, we get
\begin{equation}\label{eq:bi-fourier-delta}
    \psi_\delta \mathscr{F}[u_R]= \mathscr{F}[\Xi_{\delta,R}] +\mathscr{F}[u\opdb{\chi_R}{\mu_{\delta,1}}] +\mathscr{F}[\chi_R\opdb{u}{\mu_{\delta,1}}].
\end{equation}
with, for $\xi \in \R$, 
$$\psi_\delta(\xi):=\int_{\R}(\cos(z\xi)-1)\dd \mu_{\delta,1}(z).$$
Observe that since for all $\xi$, the map $z\mapsto \cos(z\xi)-1$ is smooth and $|\cos(z\xi)-1|\lesssim \xi^2 z^2$, by Lebesgue's theorem  $\psi_\delta$ tends to $\psi$ as $\delta\to 0$. In addition, since $u_R\in \sf H^2(\R)$, we have $\psi_{\delta}\mathscr{F}[u_R]\to \psi_{0,1}\mathscr{F}[u_R]$ in $\sfL^2(\R)$.
On the other hand, $\mathscr{F}[u\opdb{\chi_R}{\mu_{\delta,1}}]$ tends to $\mathscr{F}[u\opdb{\chi_R}{\mu_{0,1}}]$ in  $\sfL^2(\R)$, since by Plancherel's equality we have 
$$ \left\|\mathscr{F}[u\opdb{\chi_R}{\mu_{\delta,1}}]- \mathscr{F}[u\opdb{\chi_R}{\mu_{0,1}}]\right\|_{\sfL^2}\lesssim\left\|u\left(\opdb{\chi_R}{\mu_{\delta,1}}-\opdb{\chi_R}{\mu_{0,1}}\right)\right\|_{\sfL^2}\lesssim \left\|\opdb{\chi_R}{\mu_{\delta,1}}-\opdb{\chi_R}{\mu_{0,1}}\right\|_{\sfL^2}\to 0 $$
Similarly,  by Plancherel's equality we have  
$$ \|\mathscr{F}[\Xi_{\delta,R}]-\mathscr{F}[\Xi_R]\|_{\sfL^2}\lesssim \| \Xi_{\delta,R}-\Xi_R \|_{\sfL^2} $$
and as above  by Cauchy-Schwarz with respect to the measure $z^2\mu_{0,1}$ and Jensen's inequality, 
\begin{align*} 
\| \Xi_{\delta,R}-\Xi_R \|_{\sfL^2}^2&=\int_{a}^b\left(\int_{|z|\le \delta}\int_{0}^{1}\int_{0}^{1}\chi'_{R}(x+tz)u'(x+sz)z^2\dd s\dd t\dd \mu_{0,1}(z)\right)^2\,\dd x\\
&\lesssim \|u'\|^2_{\sfL^2}\int_{|z|\le \delta}|z|^2\dd \mu_{0,1}(z)\to 0 \quad\text{ as } \quad\delta\to 0.
\end{align*}
Last, let us observe that for all $\varphi\in \mathscr{S}(\R)$, 
$$\lim_{\delta\to 0}\int_{\R}\varphi(\xi)\mathscr{F}(\chi_R\opdb{u}{\mu_{\delta,1}})(\xi)\dd \xi= \lim_{\delta\to 0} \int_{\R}\mathscr{F}(\varphi)(x)\chi_R(x)\opdb{u}{\mu_{\delta,1}}(x)\,\dd x.$$
Since  $\mathscr{F}(\varphi)\in \mathscr{S}(\R)$ and $\chi_R\in \mathscr{C}_{c}^{\infty}(\R)$, by \Cref{lem:u-unif-int}, $ \mathscr{F}(\varphi)\chi_{R} \omega[u] \in \sfL^1(\R^2,\dd x\otimes\dd\mu_{0,1})$  so by Lebesgue's theorem, we deduce 
$$ \lim_{\delta\to 0} \int_{\R}\mathscr{F}(\varphi)(x)\chi_R(x)\opdb{u}{\mu_{\delta,1}}(x)\,\dd x=\int_{\R}\mathscr{F}(\varphi)(x)\chi_R(x)\opdb{u}{\mu_{0,1}}(x)\,\dd x.$$
Therefore, 
$$\lim_{\delta\to 0}\int_{\R}\varphi(\xi)\mathscr{F}(\chi_R\opdb{u}{\mu_{\delta,1}})(\xi)\dd \xi=\int_{\R}\mathscr{F}(\varphi)(x)\chi_R(x)\opdb{u}{\mu_{0,1}}(x)\,\dd x=\int_{\R}\varphi(\xi)\mathscr{F}(\chi_R\opdb{u}{\mu_{0,1}})(\xi)\dd \xi.$$

Hence, passing to the limit in \eqref{eq:bi-fourier-delta}
in $\mathscr S'(\R)$, we obtain
\[
\psi_{0,1}\mathscr{F}[u_R]
=
\mathscr{F}[\Xi_R]
+\mathscr{F}\left[u\opdb{\chi_R}{\mu_{0,1}}\right]
+\mathscr{F}\left[\chi_R\opdb{u}{\mu_{0,1}}\right]
\qquad\text{in }\mathscr S'(\R).
\]
Moreover, all the terms in the above identity being locally integrable, the equality holds almost everywhere in $\R$.
 \end{proof}
With this lemma at hand, since $\chi_R\opdb{u}{\mu_{0,1}}$ and $\Xi_R$ are in $\sf H^1(\R)$,  $\mathscr{F}[\chi_R\opdb{u}{\mu_{0,1}}]$ and $\mathscr{F}[\Xi_R]$ are in $ \mathsf{L}^1(\R)$ and thus  $\psi_{0,1} \mathscr{F}[u_R] \in \sfL^1(\R)$.  Hence, $u_R\in \mathcal{R}_{\psi_{0,1}}$ and 
$\opdb{u_R}{\mu_{0,1}}\in \mathscr{C}(\R)$. In particular, by using the definition of $\chi_R$ this implies that
$\opdb{u}{\mu_{0,1}}\in \mathscr{C}((R-1,R+1))$. $R$ being arbitrary, this implies that 
$\opdb{u}{\mu_{0,1}} \in \mathscr{C}(\R)$.
Hence we have 
\begin{equation}\label{eq:bistable-u}
    \opdb{u}{\mu} +  c \, u' + f(u) = 0,
\end{equation}
  everywhere.
\end{proof}

To conclude our construction,  we need to show that $u$ has the right limits at $\pm\infty$.

\begin{lemma}\label{lem:ubc}
    The profile $u$ constructed in \Cref{lem:convbi} satisfies, 
\begin{equation*}
    \lim_{x \to - \infty} \, u(x) = 1, \qquad \lim_{x \to + \infty} \, u(x) = 0.
\end{equation*}    
\end{lemma}

\begin{proof}[{\bf Proof of \Cref{lem:ubc}}]

Since  $u$ is bounded and nonincreasing, $\lim_{x\to \pm \infty}u(x) :=l^{\pm}$ exist and  since $u(0)=\frac{\theta}{2}$, one has $l^+\le \frac{\theta}{2}\le l^-$. 
We first show that 
\begin{equation} \label{eq:bi-lim-esti-lpm}
l^{\pm}\in \{0,\theta,1\}.
\end{equation}

To show this, since $f$ satisfies \eqref{eq:bistable}, it is sufficient to prove that $f(l^{\pm})=0$. 
Let $\rho_\tau$ be a positive mollifier to be made precise later on,  since $u$ satisfies \eqref{eq:bistable-u-pp}, we have 
\begin{equation}\label{eq:bistable-psi-tau}
c\psi'_\tau + \opdb{\psi_\tau}{\mu}+\rho_\tau\star f(u)=0,
\end{equation}
where  $\psi_\tau:=\rho_\tau \star u$.  
Note that since $u'\in \mathsf{L}^{\infty}(\R)$ and $u$ is bounded non increasing, $\psi_\tau$ is non increasing and bounded and we have for all $\tau>0$, 
\[
 \lim_{x \to \pm \infty}\psi'_\tau(x)=0, \qquad\text{ and }\qquad \lim_{x\to\pm \infty}\opdb{\psi_\tau}{\mu}(x)=0.
\]
As a consequence, by taking the limits $x\to\pm \infty$ in \eqref{eq:bistable-psi-tau}, since $u\in \mathscr{C}(\R)$, we get 
$$0=\lim_{x\to\pm\infty}\rho_\tau\star f(u(x))=\lim_{x\to\pm\infty}\int_{\R}\rho_\tau(-z)f(u(x+z))\,\dd z =f(l^{\pm}).$$

By \eqref{eq:bi-lim-esti-lpm}, we deduce readily that $l^+=0$ since  $l^+\le \frac{\theta}{2}$. 
To conclude our proof, it remains to show that $l^{-}\neq\theta$. We argue by contradiction, assume that $l^-=\theta$, then $u$ solves 
\begin{align*}
    &\opdb{u}{\mu} +  cu' + f(u) = 0.\\
   &\lim_{x\to -\infty} u(x)=\theta, \qquad\lim_{x\to +\infty} u(x)=0 
\end{align*}    
Since $u'\in \sfL^2(\R)$ and $\opdb{u}{\mu} +  cu' + f(u) \in \sfL^{\infty}(\R)$, by  multiplying \eqref{eq:bistable-u-pp} by $u'$ and integrating over $\R$, we get 
\begin{align*}
 &c\int_{\R}(u')^{2}+\int_{\R}u'\opdb{u}{\mu}(x)\,\dd x+\int_{\R}u'(x)f(u(x))\,\dd x= 0 \\
&c\int_{\R}(u')^{2} +\int_{\R}u'\opdb{u}{\mu}(x)\,\dd x+\int_{\theta}^{0}f(s)\,\dd s= 0.
\end{align*}
Let us estimate the last integral. By decomposing $\mu =\mu_{\delta}+\mu_{0,\delta}$, we get 
$$ \int_{\R}u'\opdb{u}{\mu}(x)\,\dd x=\int_{\R}u'\opdb{u}{\mu_\delta}(x)\,\dd x +\int_{\R}u'\opdb{u}{\mu_{0,\delta}}(x)\,\dd x. $$
Since $\mu_\delta$ is a finite Borel measure, by reproducing the argument in \Cref{prop:esti}, we deduce that 
$$ \int_{\R}u'\opdb{u}{\mu_\delta}(x)\,\dd x =0,$$
and therefore we deduce that  
\begin{align*}
\left|\int_{\R}u'\opdb{u}{\mu}(x)\,\dd x\right|= \left|\int_{\R}u'\opdb{u}{\mu_{0,\delta}}(x)\,\dd x\right|&= \left|\int_{\R}u'\int_{0}^1\int_{0}^1\int_{\R}u''(x+tzs)tz^2\dd \mu_{0,\delta}(z)\dd s\dd t\dd x\right|\\
&\le \frac{1}{2} \|u'\|_{\sfL^2}\|u''\|_{\sfL^2}\int_{|z|\le\delta}z^2\dd\mu(z). 
\end{align*}
Since $\mu$ is a Lévy measure, we can then find $\delta$ small such that 
$$ \left|\int_{\R}u'\opdb{u}{\mu}(x)\,\dd x\right|\le \frac{c\|u'\|_{\sfL^2}^2}{2}.$$
Thus we get 
$$\frac{c\|u'\|^2_{\sfL^2}}{2} =-\left(\int_{\R}u'\opdb{u}{\mu}(x)\,\dd x+\frac{c\|u'\|^2_{\sfL^2}}{2}\right)+\int^{\theta}_{0}f(s)\,\dd s\le  0,$$
 which contradicts the estimate $c\ge \underline{c}>0$. Hence, $l^{-}=1$.

\end{proof}

\subsection{Uniqueness of the speed} 

In this last section, we complete the proof of \Cref{thm:bi} by proving the uniqueness of the speed $c$.
When $\mu$ is a finite Borel measure, the uniqueness has already been proved in \cite{Yagisita2009a,Coville2006}, so we only need to consider uniqueness for unbounded measure $\mu$. For the rest of the section, $\mu$ is an unbounded measure. 
Assume by contradiction that there exist two solutions $(c_1,u_1)$ and $(c_2,u_2)$ to \eqref{eq:TW}. Without any loss of generality, we may assume $c_1 >c_2 >0$. Thus,  
\begin{align*}
&\opdb{u_1}{\mu} + c_1 u_1' +f(u_1)=0,\\
&\opdb{u_2}{\mu} + c_1 u_2' +f(u_2)=(c_1-c_2)u_2'\le 0, \\
&\lim_{x \to \pm \infty} (u_1-u_2)(x)=0.
\end{align*}
Since $f$ is of type \eqref{eq:bistable}, by \Cref{sec:tool}, \Cref{thm:nonlin-cp}, there exists $h_0 \in \R$ such that for all $h\le h_0$, 
$u_2(\cdot+h)\ge u_1(\cdot)$.
Define  
$$
\bar{h}:=\sup\{h \,|\,u_2(\cdot+ \tilde h)\ge u_1(\cdot), \; \text{ for all }\, \tilde h\le h  \}.
$$
We get immediately a contradiction  if  $\bar h=+\infty$, so $\bar{h}<+\infty$.
Consider $w(\cdot):=u_2(\cdot+\bar{h}) - u_1(\cdot) $, then $w\ge0$  satisfies
\begin{align*}
&\opdb{w}{\mu} + c_1 \, w' -\sup_{s\in [0,3]}|f'(s)| w \le f(u_1)-f(u_1+w) -\sup_{s\in [0,3]}|f'(s)| w\le 0  ,\\
&\lim_{x \to \pm\infty} w(x)= 0.
\end{align*}

If $w$ vanishes at $x_0\in \R$, then $x_0$ is a minimum point of $w$ (since $w \geq 0$ by construction). As a consequence, $w'(x_0) = 0$ and, furthermore by the above equation  $\opdb{w}{\mu}(x_0) \le  0$. \Cref{prop:pre-max} implies that 
$w\equiv 0$ and we have 
$$
0=\opdb{u_1}{\mu}(\cdot)+c_1 u_1'(\cdot) +f(u_1(\cdot))=\opdb{u_2}{\mu}(\cdot)+c_1 u_2'(\cdot) +f(u_2(\cdot))=(c_1-c_2)u_2'(\cdot)\le 0.
$$
Thus we have $u_2'\equiv 0 $ which contradicts  $u_2(-\infty)>u_{2}(+\infty)$.
The function $w$ is then positive, which means that $u_2(\cdot+\bar{h})>u_1(\cdot)$ so that \Cref{lem:cp} can be used to obtain a contradiction by constructing $h^{\#} > \bar{h}$  such that $u_2(\cdot+\tilde{h}) \geq u_1(\cdot)$ for all $\tilde h\le h^{\#}$, contradicting the definition of $\bar{h}$.
As a consequence, the speed is unique.
To obtain the uniqueness of the profile,  we use a similar argument and assume that two profiles $u_1$ and $u_2$ exist. As above, by  \Cref{thm:nonlin-cp}, there exists $h_0 \in \R$ such that for all $h\le h_0$, 
$u_2(\cdot+h)\ge u_1(\cdot)$ and we can define  
$$
\bar{h}:=\sup\{h \,|\,u_2(\cdot+ \tilde h)\ge u_1(\cdot), \; \text{ for all }\, \tilde h\le h  \}.
$$
 Again $\bar{h}<+\infty$ since otherwise we get a contradiction. As above consider $w(\cdot):=u_2(\cdot+\bar{h}) - u_1(\cdot) $, then $w\ge0$  satisfies
\begin{align*}
&\opdb{w}{\mu} + c \, w'  = f(u_1) -f(u_1+w),\\
&\lim_{x \to \pm\infty} w(x)= 0.
\end{align*}
Observe that by definition $w\in\mathscr{C}(\R)$ and either $w>0$ in $\R$ or  $w$ vanishes at $x_0\in \R$. In the first situation, then 
$u_2(\cdot+\bar{h})>u_1(\cdot)$ so that \Cref{lem:cp} can be used to obtain a contradiction by constructing $h^{\#} > \bar{h}$  such that $u_2(\cdot+h) \geq u_1(\cdot)$ for all $h\le h^{\#}$, contradicting the definition of $\bar{h}$. Therefore, this situation is ruled out and necessarily there exists a point $x_0$ at which $w$ attains its minimum (since $w \geq 0$ by construction). As a consequence, $w'(x_0) = 0$ and, furthermore, $\opdb{w}{\mu}(x_0) \leq 0$. By using recursively \Cref{prop:pre-max} we then  deduce that $w\equiv 0$, showing thus that $u_2(\cdot+\bar{h})\equiv u_1(\cdot)$. 

\section{Ignition type nonlinearities - Proof of \texorpdfstring{\Cref{thm:ign}}{thm:ign}}\label{sec:ign}

As in the bistable situation, the only remaining case to be considered is the one where $\mu$ is unbounded and satisfies the moment condition.

The construction in this situation essentially follows the arguments developed in \Cref{sec:bi}. For $\eps>0$, the existence of a solution $(c_\eps,u_\eps)$ to \eqref{eq:ueps} with a positive speed $c_\eps$ is then guaranteed by the results in \cite{Coville2007d}. 
To extract from $(c_\eps,u_\eps)_{\eps\in (0,1)}$ a converging subsequence, we first obtain some uniform bounds on $(c_\eps)_{\eps > 0}$. Start with an estimate from below. 
\begin{prop}\label{prop:esti-ig1}
Assume that $f$ is of type \eqref{eq:ignition} and $\mu$ is unbounded such that $\ds{\int_{\R}\min(|z|,z^2)\,\dd \mu(z)<+\infty}$. Then there exists $\eps_0>0$ and $\underline{c}>0$ such that for all $0<\eps<\eps_0$,  $$c_\eps\ge \underline{c}.$$ 
\end{prop}

\begin{proof}[{\bf Proof of \Cref{prop:esti-ig1}}]
Take a function $f_b\in \mathscr{C}^1(\R)$ of type \eqref{eq:bistable} such that $f\ge f_b$ and $\int_{0}^1f_b(s)\,ds>0$. \Cref{thm:bi} applied to $f_b$ leads to a solution with a $c_{b,\eps}$ uniformly bounded from below (by \Cref{lem:esti3}) for $\eps\in(0,\eps_0)$ for some $\eps_0>0$. From the monotonicity of $c_\eps$ with respect to the nonlinearity \cite{Coville2005}, $c_\eps\ge c_{b,\eps} \geq \underline{c}$.
\end{proof}

Next, we obtain an upper bound  for $c_\eps$. Namely, we show 

\begin{prop}\label{prop:esti-ig2}
Assume that $f$ is of type \eqref{eq:ignition} and $\mu$ is unbounded such that $\ds{\int_{\R}\min(|z|,z^2)\,\dd \mu(z)<+\infty}$.  For $0<\eps<\eps_0$, let  $(c_\eps,u_\eps)$ be  a  solution to \eqref{eq:ueps} normalised by $u_{\eps}(0)=\theta$. Then  $c_\eps$ satisfies
 $$  c_{\eps} \lesssim \int_\R \min(\vert z \vert,z^2) \, \dd \mu(z).$$
\end{prop}   

\begin{proof}[{\bf Proof of \Cref{prop:esti-ig2}}]
Combining \Cref{prop:esti} -- \eqref{1}-\eqref{4} and \Cref{prop:esti-ig1}, we have a uniform $\mathsf{H}^1(\R)$  bound for $(u'_\eps)_{\eps > 0}$. As a consequence, by Sobolev embedding, the family $(u'_\eps)_{\eps > 0}$ is also uniformly bounded in $\mathsf{L}^{\infty}(\R)$.
To estimate $\opdb{u_\eps}{\mu_\eps}$, write
\begin{align*}
   \opdb{u_\eps}{\mu_\eps} (\cdot)&=  \int_{|z|\ge \eps}z \int_{0}^1u_\eps'(\cdot+sz)\dd s\,\dd \mu(z)\\
   &= \int_{1 \geq |z|\ge \eps}z \int_{0}^1 u_\eps'(\cdot+sz)\dd s\,\dd \mu(z) +  \int_{|z|\ge 1} z \int_{0}^1u_\eps'(\cdot+sz)\dd s\,\dd \mu(z)\\
   &= \int_{1 \geq |z|\ge \eps}z \left( \int_{0}^1 \left( u_\eps'(\cdot) + sz\int_{0}^{1} u_\eps''(\cdot+s\nu z) \, \dd \nu \right) \dd s \right) \,\dd \mu(z) +  \int_{|z|\ge 1} z \int_{0}^1u_\eps'(\cdot+sz)\dd s\,\dd \mu(z)\\
   &= \int_{1 \geq |z|\ge \eps} \int_{0}^1 \int_{0}^{1} sz^2 \, u_\eps''(\cdot+s\nu z) \, \dd \nu \, \dd s \,\dd \mu(z) +  \int_{|z|\ge 1} z \int_{0}^1u_\eps'(\cdot+sz)\dd s\,\dd \mu(z),
\end{align*}
where we have used twice Taylor's theorem with $u_\eps\in \mathscr{C}^2(\R)$ and the symmetry of $\mu$. Observe that Fubini's theorem combined with the fact that $\lim_{x \to + \infty} u_\eps(x) = \lim_{x \to + \infty} u_{\eps}'(x)=0$ gives
\begin{align}
    \int_{\R_+} \opdb{u_\eps}{\mu_\eps}(x) \, \dd x &= - \int_{1 \geq |z|\ge \eps} \int_{0}^1 \int_{0}^{1} sz^2 \, u_\eps'(s\nu z) \, \dd \nu \, \dd s \,\dd \mu(z) -  \int_{|z|\ge 1} z \int_{0}^1u_\eps(sz)\dd s\,\dd \mu(z)\label{eq:theta-ceps-form1}\\
    &\lesssim \int_{\R} \min(\vert z \vert , z^2) \dd \mu(z)\nonumber.
\end{align}
Since $u_\eps(0)=\theta$ and $u_\eps$ is monotone decreasing  we get $f(u_\eps)= 0$ on $\R_{+}$. Integrate \eqref{eq:ueps} over $\R_{+}$, to get 
\begin{equation}\label{eq:theta-ceps-form2}
 c_\eps \theta =\int_{\R_{+}} \opdb{u_\eps}{\mu_\eps}(y) \, \dd y \lesssim \int_{\R} \min(\vert z \vert , z^2) \dd \mu(z). 
\end{equation}
\end{proof}
With these estimates at hand, by normalising $u_\eps$ by $u_\eps(0)=\theta$, we may reproduce all the strategy explained in \Cref{sec:bi}. From this, there exists $(c,u) \in \R_{+}^{\star} \times \mathscr{C}^{1,\alpha}(\R)$ for all $\alpha \in \left(0,\frac{1}{2}\right)$.

The last difficulty is to prove that $u$ satisfies the right limits at $\pm\infty$, since the uniqueness of the speed in this situation is obtained as in the bistable case. Again, since $u$ is bounded nonincreasing, $u$ achieves some limits at $\pm \infty$ called $l^{\pm}$ and, arguing as in the bistable case, these limits satisfy $l^+\le \theta\le l^-$ with $f(l^{\pm})=0$. 

We first show that $l^{-}=1$. Necessarily, since $u$ is nonincreasing, $u(0)=\theta$, and $f$ is of type \eqref{eq:ignition}, $f(l^{-})=0$ implies that $l^{-}\in \{\theta,1\}$. We now argue by contradiction and assume that $l^-=\theta$. Then $u$ satisfies
$$
cu'+ \opdb{u}{\mu}= 0,
$$
on $\R$ and by monotonicity of $u$, $u\le \theta$ and $u$ reaches a maximum at $0$. By \Cref{prop:pre-max}, this implies that $u\equiv \theta$. Therefore, the sequence $(u_{\eps})_{\eps > 0}$ converges to $\theta$ pointwise and $(u_{\eps}')_{\eps > 0}$ tends to $0$ locally uniformly.

Since, by \eqref{eq:theta-ceps-form1} and \eqref{eq:theta-ceps-form2} we have
\begin{equation*}
\theta c_\eps=-\int_{|z|\le 1}\iint_{0}^1sz^2u'_\eps(\nu sz)\, \dd s\dd \nu \dd \mu_{\eps}(z) -\int_{|z|\ge 1}\int_{0}^1zu_\eps(sz)\, \dd s\dd \mu_{\eps}(z),
\end{equation*}

sending $\eps\to 0$ and using that $u_\eps \to \theta$ pointwise and  $u'_\eps \to 0$ locally uniformly, we then obtain the following contradiction 
$$
0<\theta c =- \theta \int_{|z|\ge 1}z\dd \mu(z)=0.
$$
Therefore  $l^-=1$. To prove that $l^{+}=0$, we actually take advantage of the previous relation to observe that, since $u_\eps \to u$ pointwise and locally uniformly in $\mathscr{C}^{1,\alpha}(\R)$ and  $c_\eps$ converges to $c$, 
\begin{equation}\label{eq:ign-lim1}
c\theta= -\iint_{0}^{1}\int_{|z|\le 1}sz^2u'(\nu sz)\, \dd s\dd \nu \dd \mu(z) -  \int_{0}^{1}\int_{|z|\ge 1}z u(sz)\, \dd s\dd \mu(z).
\end{equation}
On the other hand, since $u\in \mathscr{C}^1(\R)$, by using the symmetry of $\mu$, $u$ satisfies
$$
-cu'(\cdot)=\int_{0}^{1}\int_{|z|\le 1}z [u'(\cdot+sz)-u'(\cdot)]\dd s\dd \mu(z) + \int_{0}^{1}\int_{|z|\ge 1}zu'(\cdot+sz)\, \dd s\dd \mu(z) +f(u(\cdot)).
$$
Integrating the above equation over $\R_{+}$, since $u(0)=\theta$, 
we get by using Fubini's theorem,
\begin{align*}
c(\theta -l^+)&=\int_{0}^{1}\int_{|z|\le 1}z [\theta - u(sz)] \, \dd s\dd \mu(z)  -\int_{0}^{1}\int_{|z|\ge 1}z u(sz)\, \dd s\dd \mu(z) +l^+\int_{0}^{1}\int_{|z|\ge 1}z \dd s\dd \mu(z),
\end{align*}
which after using the Taylor theorem in the first integral and the symmetry of $\mu$  reduces to 
\begin{equation*}
c(\theta -l^+)= -\iint_{0}^{1}\int_{|z|\le 1}sz^2u'(\nu sz)\, \dd s\dd \nu \dd \mu(z) -  \int_{0}^{1}\int_{|z|\ge 1}z u(sz)\, \dd s\dd \mu(z) = c \theta,
\end{equation*}
recalling \eqref{eq:ign-lim1}. We deduce $-c l^+ =0$, which implies that  $l^+=0$.

\section{Monostable type nonlinearities - Proof of \texorpdfstring{\Cref{thm:mono}}{thm:mono}} \label{sec:mono}
In this last section, we prove \Cref{thm:mono}. To do so, let us consider a smooth monotone cut-off function $\chi\in \mathscr{C}^{\infty}(\R)$ such that
$$
\chi(x):=\begin{dcases}
1 \quad \text{ for } \quad x\ge \frac{1}{2}\\
0 \quad \text{ for } \quad x\le \frac{1}{4}.
\end{dcases}
$$ 
For $n \in \N^\star$, consider $\chi_n(\cdot):=\chi\left(n \,\cdot\right)$ and define the following $\mathscr{C}^{1}$ function $f_n:=f \, \chi_n$.
By definition $f_n$ is of type \eqref{eq:ignition} and from \Cref{thm:ign}, there exists $(c_n,u_n) \in \R_{+}^{\star} \times \mathscr{C}^1(\R)$, solution to 
\begin{align*}
&\opdb{u_n}{\mu} +c_n u_n' +f_{n}(u_{n})=0,\\
&\lim_{x\to -\infty} u_{n}(x)=1 \qquad \lim_{x\to +\infty} u_{n}(x)=0.
\end{align*}

Observe that, by definition of $f_n$, we have $f_n\le f_{n+1}\le f$ for all $n\in \N$. Using the definition of $(\kappa,w)$, we also have, for all $n >0$,
\begin{align*}
&\opdb{w}{\mu} +\kappa w' +f_{n}(w)\le  \opdb{w}{\mu} +\kappa w' +f(w)\le 0.\\
&\lim_{x \to -\infty} w (x) \geq 1 = \lim_{x \to -\infty} u_{n}(x),\\
&\lim_{x \to +\infty} w (x) \geq 0 = \lim_{x \to +\infty} u_{n}(x).
\end{align*}
As a consequence, from \Cref{thm:mono-c} we deduce that   $(c_n)_{n\in\N}$ is an increasing sequence and that $c_n\le \kappa$ for all $n \in \N^\star$. Consequently, $(c_n)_{n \in \N}$ converges to some $0< c^\star \le \kappa $. Arguing as in \Cref{sec:bi}, after normalising $u_n(0)=\frac{1}{2}$, we can pass to the limit $n$ to $+\infty$ to obtain a  monotone bounded   solution $(c^\star,u) \in \R_{+}^{\star} \times \mathscr{C}^1(\R)$ to  
$$
\opdb{u}{\mu}+c^\star u'+f(u)=0.
$$
The limits of $u$ at $\pm \infty$ immediately follow from the fact that $f(l^{\pm})=0$, that can be obtained as in the proof of \Cref{lem:ubc}, and $l^-\ge \frac{1}{2}\ge l^{+}$.

To construct front solutions for $c>c^\star$, we introduce the following problem
\begin{align*}
&\opdb{v}{\mu_\eps} +c_{n,\eps} v' +f_{n}(v)=0, \qquad \text{ on }\; \R, \\
&\lim_{x\to -\infty} v(x)=1 \qquad \lim_{x\to +\infty} v(x)=0 
\end{align*}

Thanks to the previous section, for fixed $n \in \N^\star$ and $\eps>0$, there exists $(c_{n,\eps},u_{n,\eps})$, solution to that problem.

For $\eps>0$, define $c^{\#}(\eps):=\lim_{n \to + \infty} c_{n,\eps}.$ 
Then either there exists a sequence $(\eps_k)_{k\in\N}$ of positive numbers such that $\eps_k\to0$ as $k\to+\infty$ and
$c^{\#}(\eps_k)=+\infty$ or there exists $\eps_0$ such that for all $\eps\le\eps_0$, $c^{\#}(\eps)<+\infty$.
Assume that the sequence $(\eps_k)_{k\in\N}$ exists, then, there exists a sequence $(c_{\varphi(k),\eps_{k}})_{k\in\N}$ such that  
$\eps_{k} < \eps_0$ and
$c_{\varphi(k),\eps_{k}}>2c$ for all $k\ge k_0$. 
By uniqueness of the speed, we have also, for all fixed $k \in \N$, 
$$\lim_{\eps \to 0}c_{\varphi(k),\eps}=c_{\varphi(k)} \leq c^\star,$$
and thus there exists $\eps_{\varphi(k)}' \in (0,\eps_k)$ such that $c_{\varphi(k),\eps_{\varphi(k)}'}\leq \tfrac{1}{2}(c+c^\star) < c$. 
By \cite{Coville2006}, we can check that the map $\eps \mapsto c_{\varphi(k),\eps}$ is continuous on $(\eps_{\varphi(k)}',\eps_{k})$, and by the intermediate value lemma, for any $k \geq k_0$, there exists $\overline{\eps}_{\varphi(k)} \in (\eps_{\varphi(k)}',\eps_{k})$ such that $c_{\varphi(k),\overline{\eps}_{\varphi(k)}} = c$. 
There is then a sequence $(\bar \eps_{\varphi(k)})_{k \in\N}$ converging to $0$ and a solution $u_{\varphi(k),\bar\eps_{\varphi(k)}}$ with $u_{\varphi(k),\bar \eps_{\varphi(k)}}(0)=\frac{1}{2}$, such that after sending $k\to +\infty$, and working as in \Cref{ss:constr-front} we get a non trivial solution to  \eqref{eq:TW} with speed $c$.

Assume now that this sequence $(\eps_k)_{k\in\N}$ does not exist. In this situation,for all $\eps\le \eps_0$, $\lim_{n\to+\infty}c_{n,\eps}=c^{\#}(\eps)<+\infty$  and there exists a positive $u_\eps$ solution to 
\begin{align*}
&\opdb{u_\eps}{\mu_\eps} +c^{\#}(\eps) u_\eps' +f(u_\eps)=0 \qquad \text{ on }\, \R,\\
&\lim_{x\to -\infty} u_\eps(x)=1, \quad \lim_{x\to +\infty} u_\eps(x)=0.
\end{align*}
From \cite{Coville2007a}, we can also infer that the above system possesses front solutions of any speed $c\ge c^{\#}(\eps)$. Therefore, for all $c\ge\liminf_{\eps\to 0}c^{\#}(\eps)$, by sending $\eps\to0$ and working as in \Cref{ss:constr-front} we get a non trivial solution to  \eqref{eq:TW} with speed $c$. 

To complete our proof we may reproduce the scheme of proof developed in \cite{Coville2007a}. 
Namely, by using the monotonicity of the speed, we have the following diagram 
\[ \begin{tikzcd}
c_n \arrow{r}{n\to\infty}  & c^{\star} \le   \liminf c^{\#}(\eps)  \\%
c_{n,\eps} \arrow{u}{\eps\to0} \arrow[swap]{r}{n\to +\infty}&  c^{\#}(\eps) \arrow{u}{\eps \to 0}
\end{tikzcd}
\]
If $c^\star= \liminf_{\eps\to0} c^{\#}(\eps)$ then there is nothing to prove so let us assume that $c^{\star} <  c^{\star\star}:=\liminf c^{\#}(\eps) $ and pick $c\in (c^{\star},c^{\star\star})$. Then from the  diagram, we can reproduce the argument used above to find a sequence $(\eps_n)_{n\in\N}$  such that $\eps_n\to 0, c_{n,\eps_n}=c$ with $u_{n,\eps_n}(0)=\frac{1}{2}$, such that after sending $n$ to $+\infty$, and working as in \Cref{ss:constr-front} we get a non trivial solution to  \eqref{eq:TW} with speed $c$.

\section{Acknowledgements}
The authors are deeply grateful to Professors Petru Mironescu and Jean-Yves Chemin for fruitful discussions on measure theory and Fourier analysis that led to substantial simplification and clarification of the manuscript. E. Bouin and J. Coville acknowledge support from the ANR via the project REACH under grant agreement ANR-23-CE40-0023-01. 

\appendix

\section{Mathematical toolbox}\label{sec:tool}

\subsection{Linear comparison principles}

To be reasonably self-contained, we recall below a few basic comparison principles for sub- and super- solutions associated to a generator $\mathcal{L}$, in the elliptic and parabolic cases. These being standard, we do not provide proofs in this document. However, an interested reader could come up with a proof by adapting some from the literature, for example from \cite{Andreu-Vaillo,Barles2008a,Coville2007,Coville2008} for the comparison principle  and  \cite{Andreu-Vaillo,Brasseur2021,Cabre2013,Chen2002,Zhang2023} for its parabolic version.   

The following statement is crucial in the proofs of comparison principles. 
\begin{prop}\label{prop:subgrp}{\cite{Alibaud2018,Choquet1960,Coville2008}}\label{prop:pre-max}
Let $u\in \mathscr{C}(\R)\cap \mathsf{L}^{\infty}(\R)\cap \text{Dom}(\D_\mu)$ be such that 
$$
\opdb{u}{\mu}(x) \ge 0  \quad \text{ for all }\quad x\in\R. 
$$
Assume that $u$ achieves a nonnegative global maximum at some point $x_0 \in\R$, then 
  for all $q\in \text{Supp}(\mu)$,  
 $$x_0+ q\,\N  \subset \{x \in \R\,|\,u(x)=u(x_0)\}.$$ 
 When $\mu$ is symmetric and unbounded, then  $\{x \in \R\,|\,u(x)=u(x_0)\}=\R$. 
\end{prop}

\begin{proof}[{\bf Proof of \Cref{prop:subgrp}}]
The proof of the first part  is a trivial consequence of the definition of being a maximum and induction.  When $\mu$ is symmetric and the domain $\R$ then it is shown in  \cite{Alibaud2018,Choquet1960} that  $\{x \in \R\,|\,u(x)=u(x_0)\}$ contains  a sub-group of $\R$ which is dense  when $\mu(\R)=+\infty$. Since $\{x \in \R\,|\,u(x)=u(x_0)\}$  is a closed set, we are done.
\end{proof}
A similar proposition can be obtained on $\R^+$. Namely,
\begin{prop}\label{prop:subgrp2}{\cite{Coville2008}}\label{prop:si-pre-max}
Let $u\in \mathscr{C}(\R)\cap \mathsf{L}^{\infty}(\R)\cap \text{Dom}(\D_\mu)$ be such that 
$$
\opdb{u}{\mu}(x) \ge 0  \quad \text{ for all }\quad  x\in \R_{+}.
$$
Assume that $u$ achieves a nonnegative global maximum at some point $ x_0\in \R_{+}$, 
then,  for all $q\in \text{Supp}(\mu)$,  
 $$(x_0+ q\,\N)\cap \R_+  \subset \{x \in \R_{+}\,|\,u(x)=u(x_0)\}.$$ 
 When $\mu$ is symmetric and unbounded, then    $\{x \in \R_+\,|\,u(x)=u(x_0)\}=\R_+$. 
\end{prop}

From \Cref{prop:pre-max} and \Cref{prop:si-pre-max}, we can derive the following comparison theorem
\begin{theorem}[Elliptic type comparison principle]\label{thm:cp}

\qquad\medskip

Let $b,d\in \mathscr{C}(\R)$ with $d\le0$.  Assume that $u,v \in \textrm{Dom}(\D_\mu)\cap\mathscr{C}^1(\R)$ are such that either
\begin{align*}
&\opdb{u}{\mu}+bu'+du\ge 0 \quad \text{ and } \quad \opdb{v}{\mu}+bv'+dv\le 0,\quad\text{ on } \quad \R,\\[5pt]
&\lim_{x\to \pm \infty} (u-v)(x) \leq 0,
\end{align*}
or,
\begin{align*}
&\opdb{u}{\mu}+bu'+du\ge 0 \quad \text{ and } \quad \opdb{v}{\mu}+bv'+dv\le 0,\quad\text{ on } \R_{+},\\[5pt]
&u\le v \quad\text{ on } \quad x\in\R^{\star}_{-} \quad \text{ and } \quad \lim_{x\to +\infty} (u-v)(x) \leq 0.
\end{align*}
Then $u \le v$ on $\R$.
\end{theorem}

\begin{theorem}[Parabolic type comparison principle]\label{thm:pcp}

\qquad\medskip

Let $u$ and $v$ be two functions of class $ \mathscr{C}^1((0,T), \textrm{Dom}(\D_\mu) )\, \cap \, \mathscr{C}([0,T],\textrm{Dom}(\D_\mu))$ such that
$u$, $v$, $\partial_t u$, $\partial_t v$ are bounded and 
\begin{align*}
&\partial_t u -\opdb{u}{\mu} \le 0 \le \partial_t v  -\opdb{v}{\mu}, &\quad\text{ on }\quad \left\lbrace (t,x) \in (0,T] \times \R, x > x(t) \right\rbrace,&\\[5pt]
&u(0,\cdot)\le_{\not \equiv} v(0,\cdot),\\[5pt]
&u(t,\cdot)< v(t,\cdot), &\quad\text{ on }\quad \left\lbrace (t,x) \in (0,T] \times \R, x \le x(t) \right\rbrace,\\[5pt]
&\lim_{\pm \infty} u(t,x)-v(t,x) \leq 0,& \quad t \in(0,T].
\end{align*}
Then $u \le v$ on $[0,T] \times \R$.
\end{theorem}

\subsection{Sliding and consequences}

From the above comparison principle, using sliding techniques in the spirit of \cite{Berestycki1991,Coville2006,Vega1993}, rigidity results can be derived for semi-linear equations of the form   
\begin{equation*}
 cu' +\opdb{u}{\mu} +f(u )=0,    
\end{equation*}
on $\R$. Namely,
\begin{theorem}\label{thm:nonlin-cp}
Let $f\in \mathscr{C}^1(\R)$ be such that for $\gamma>0$, $f(0)=f(\gamma)=0$ and $f<0$ in $(\gamma,+\infty)$.  Assume that there exists $\delta_0>0$ such that $f$ is non increasing  in $s\in[0,\delta_0]\cup [\gamma-\delta_0,\gamma]$.  Let $\overline{u}$ and $\underline{v}$ be two functions of class $\mathsf{L}^{\infty}(\R)\cap \mathscr{C}^{1}(\R,\R_+)\cap \textrm{Dom}(\D_\mu)$, (of class $\mathsf{L}^{\infty}(\R)\cap \mathscr{C}(\R,\R_+)\cap \textrm{Dom}(\D_\mu)$ when $c=0$),  such that 
\begin{align*}
&\D_\mu[\overline{u}]+c\overline{u}' +f(\overline{u})\le 0 \leq \D_\mu[\underline{v}]+c\underline{v}' +f(\underline{v}),\\[5pt]
&\lim_{x \to +\infty} \overline{u}(x)\ge 0, \quad \lim_{x \to +\infty} \underline{v}(x)\le 0,\\[5pt]
&\lim_{x \to -\infty} \overline{u}(x)\ge \gamma, \quad \lim_{x \to -\infty} \underline{v}(x)\le \gamma.
\end{align*} 
Then there exists $\tau_0\in\R\cup \{+\infty\}$ such that for every $ \tau\in (-\infty,\tau_0) $, $u(\cdot+\tau)> v(\cdot)$. 
When $\ds{\lim_{x \to -\infty} \underline{v}(x)<\gamma}$, then $\tau_0=+\infty$.
\end{theorem}

A straightforward corollary of this nonlinear comparison principle is the monotonicity of the speed $c$ with respect to the nonlinearity. 
\begin{theorem}\label{thm:mono-c}
Let $f,g\in \mathscr{C}^1(\R)$ be such that $f\le g, f\not\equiv g$ and for $0<\gamma_f\le\gamma_g\le 1$, $f(0)=f(\gamma_f)= g(0)=g(\gamma_g)=0$. Assume that either $f$ or $g$ satisfies the assumptions of \Cref{thm:nonlin-cp}.  Let $\overline{u}$ and $\underline{v}$ in $\mathsf{L}^{\infty}(\R)\cap \mathscr{C}^{1}(\R,\R_{+})\cap \textrm{Dom}(\D_\mu)$ be non constant and let  $c_1$ and $c_2$ be  such that 
\begin{align*}
&\D_\mu[\overline{u}]+c_1 \overline{u}' +g(\overline{u})\le 0 \le \D_\mu[\underline{v}]+c_2\underline{v}' +f(\underline{v}),\\[5pt]
&\lim_{x \to +\infty} \overline{u}(x)\ge 0, \quad \lim_{x \to +\infty} \underline{v}(x)\le 0,\\[5pt]
&\lim_{x \to -\infty} \overline{u}(x)\ge \gamma_g, \quad \lim_{x \to -\infty} \underline{v}(x)\le \gamma_f.
\end{align*} 
Then, $c_1 \ge c_2$. If $g_{|(0,\gamma_g)}\not\equiv f$, then $c_1>c_2$.
\end{theorem}
 At the core of the proof of  \Cref{thm:nonlin-cp}, there is  the following technical sliding lemma:
 \begin{lemma}\label{lem:cp}
Let $f$ be as in \Cref{thm:nonlin-cp}. Let $\overline{u}$ and $\underline{v}$ be as in \Cref{thm:nonlin-cp} 
 and such that the conditions 
\begin{align*}
\gamma-\overline{u} <\tfrac{\delta_0}{2} \quad \text{ on } \quad (-\infty,-M) \quad \text{ and } \quad v(x)< \tfrac {\delta_0}{2} \quad \text{ on } \quad(M,+\infty),
\end{align*}
are satisfied for some $M>0$. If there exists a positive constant $b$ such that $\overline{u}$ and $\underline{v}$ satisfy:
\begin{align*}
\overline{u}(\cdot+b)> \underline{v} \quad \text{ on } \quad [-M-1,M+1] \quad \text{ and } \quad \overline{u}(\cdot+b)+\tfrac{\delta_0}{2}>\underline{v} \quad \text{ on } \quad \R,
\end{align*}
then we have $\overline{u}(\cdot+b)\ge \underline{v}$ everywhere.
\end{lemma}

\subsection{A technical lemma and some properties of $\textrm{Dom}(\D_{\mu})$}
We first obtain a technical lemma (a commutator estimate).

\begin{lemma}\label{lem:tec-f}
    Let $f \in \mathscr{C}^1(\R)$ be such that $f(0)=0$. Let  $u\in \mathscr{C}(\R,\R^\star_+)\cap \sfL^{\infty}(\R)$ be a nonincreasing function. Let $\rho\in\mathscr{C}^{\infty}_c(\R)$ be a positive function of unit mass with support in $(-1,1)$. For $\tau>0$, consider $\rho_\tau=\frac{1}{\tau}\rho\left(\frac{\cdot}{\tau}\right)$. Then for any $\eps>0$, there exists $\tau(\eps)$ such  that for all $\tau\le \tau(\eps)$,
    $$|\rho_\tau \star f(u)- f(\rho_\tau\star u )|\le  \eps \rho_\tau\star u.$$
\end{lemma}
\begin{proof}[{\bf Proof of \Cref{lem:tec-f}}]
Observe that since $u$ is bounded, monotone and continuous and  $f\in \mathscr{C}^1(\R)$, 
we have $|f'(u)|\le \|f'\|_{\sfL^{\infty}(0,\|u\|_{\infty})}=:K$ and
$$\lim_{\tau \to 0}\|\rho_\tau\star f(u) -f(u)\|_{\mathsf{L}^{\infty}(\R)}=\lim_{\tau \to 0}\|\rho_\tau\star u -u\|_{\mathsf{L}^{\infty}(\R)}=0.$$
From this observation, on the set  $\{u\ge \delta\}$,  since $f\in \mathscr{C}^1(\R)$ and $u\in \mathsf{L}^{\infty}(\R)$,  we have for $\tau$ small, say $\tau\le \tau_\delta$, $\rho_\tau\star u(x)\ge \frac{\delta}{2}$  and 
\begin{align*}
| \rho_\tau \star f(u)- f(\rho_\tau\star u )| &\le | \rho_\tau \star f(u)-f(u)|+ |f(u)- f(u+(\rho_\tau\star u-u) )| \\
&\le \frac{\eps\delta}{4}+K\|\rho_\tau\star u-u\|_{\mathsf{L}^{\infty}(\R)} \le \frac{\eps\delta}{2} \le \eps\rho_\tau\star u.    
\end{align*}
Let  $\ds{l^+=\lim_{x\to+\infty} u(x)}$. If $l^+>0$, then we are done by taking $\tau_\delta$ with $\delta=l^+$. 
Let us then assume that $l^+=0$. 
Since $f\in \mathscr{C}^1(\R)$ by Taylor's formula, $f(s)=f'(0)s+r(s)$ with $\lim_{s\to 0}\left|\frac{r(s)}{s}\right|=0$. Fix  $\delta_0>0$ such that $\left|\frac{r(s)}{s}\right|\le \frac{\eps}{2}$ for $s\le \delta_0$.
Since $l^+=0$ and $u$ non increasing and continuous, there exists  $r_0<r_1$ such that  $u(r_0)= \frac{\delta_0}{2}$ and $u(r_1)= \frac{\delta_0}{4}$.

Let $\tau_0\le r_1-r_0$, then for all $\tau\le \tau_0$, $\rho_\tau\star u \le \frac{\delta_0}{2}$ and 
for $x\ge r_1$,
\begin{align*} 
| \rho_\tau \star f(u)- f(\rho_\tau\star u )|= | \rho_\tau \star r(u)- r(\rho_\tau\star u )|&\le  \rho_\tau \star \left(u \left| \frac{r(u)}{u}\right|\right)+  \left|\frac{r(\rho_\tau\star u)}{\rho_\tau\star u}\right| \rho_\tau\star u\\
&\le \left(\frac{\eps}{2}  +\frac{\eps}{2}\right)\rho_\tau\star u
\end{align*}
The above estimate holds as well for $x\le r_1$, by choosing $\tau(\eps)\le \min\{\tau_0,\tau_{\delta_0/4}\}$. Hence, we have 
on $\R$
\[
| \rho_\tau \star f(u)- f(\rho_\tau\star u )|
\le \eps\rho_\tau\star u.
\]
\end{proof}

Last, we obtain a Jacob-Schilling type characterisation of $\text{Dom}(\D_\mu)$ 
\begin{lemma}\label{lem:dom-mu}
Let $\mu$ be a symmetric L\'evy measure and denote $$\ds{\psi : \xi \mapsto \int_{\R}(\cos(z\xi)-1)\dd \mu(z)}$$ its symbol. Let $\mathcal{R}_{\psi}$ be the space of functions 
$$ \mathcal{R}_{\psi}:=\{ u\in \mathscr{C}_0(\R)\cap \sfL^1(\R), \psi\mathscr{F}[u] \in \sfL^1(\R)\}.$$
Then $\mathcal{R}_{\psi}\subset \mathscr{A}:=\{u\in \mathscr{C}(\R)\cap \text{Dom}(\D_{\mu})\, |\, \opdb{u}{\mu} \in \mathscr{C}(\R)\}$, and for $u\in \mathcal{R}_{\psi}$,
$$\opdb{u}{\mu}=\mathscr{F}^{-1}\left(\psi(\xi)\mathscr{F}[u]\right). $$
\end{lemma}
\begin{proof}[{\bf Proof of \Cref{lem:dom-mu}}]
For $\delta>0$ and $u\in\mathcal{R}_{\psi}$, since $\mu_\delta$ is a finite Borel measure, we have 
$$\mathscr{F}[\opdb{u}{\mu_\delta}]= \psi_{\delta} \mathscr{F}[u], $$
with $\ds{\psi_\delta=\int_{\R}(\cos(z\xi)-1)\dd\mu_\delta(z)}$.
Since for all $\xi$, the map $z\mapsto \cos(z\xi)-1$ is smooth,  we have  $\lim_{\delta\to 0} \psi_\delta =\psi $ at least pointwise.
Besides, by definition of $\psi_\delta$, we get  $|\psi_{\delta}|\le|\psi|$ and $|\psi_\delta \mathscr{F}[ u]|\le |\psi\mathscr{F}[u]|\in \sfL^1(\R)$. By Lebesgue's theorem, it follows that  $\lim_{\delta \to 0}\|(\psi_\delta-\psi)\mathscr{F}[ u]\|_{\sfL^1(\R)} = 0$.

On the other hand, we have,
\begin{align*}
\left|\opdb{u}{\mu_\delta}-\mathscr{F}^{-1}(\psi\mathscr{F}[ u])\right| = \left|\mathscr{F}^{-1}(\psi_\delta\mathscr{F}[ u])-\mathscr{F}^{-1}(\psi\mathscr{F}[u])\right|&= \left|\mathscr{F}^{-1}([\psi_\delta-\psi]\mathscr{F}[ u])\right| \\
&\lesssim \|(\psi_\delta-\psi)\mathscr{F}[u]\|_{\sfL^1(\R)}.
\end{align*}
As a consequence, for all $x\in\R$,
$$ \lim_{\delta\to 0} \left|\opdb{u}{\mu_\delta}(x)-\mathscr{F}^{-1}(\psi\mathscr{F}[u])(x)\right|=0,$$
and thus $u\in\text{Dom}(\D_\mu)$ with 
$$\opdb{u}{\mu}=\mathscr{F}^{-1}\left(\psi(\xi)\mathscr{F}[u]\right). $$
Since the map $x\mapsto \mathscr{F}^{-1}(\psi\mathscr{F}[u])$ is continuous, $u\in \mathscr{A}$.
\end{proof}

\bibliographystyle{plain}
		\bibliography{gentw}

\end{document}